\documentclass[12pt,twoside]{article}
\usepackage{graphicx,amsmath,amsfonts,amssymb,amscd,bezier,amstext,makeidx,esint}
\usepackage{amsthm}
\usepackage{fancyhdr}
\usepackage[bottom=2cm,top=3cm,left=2cm,right=2cm]{geometry}
\usepackage[Sonny]{fncychap}
\usepackage[T1]{fontenc}
\usepackage{pst-all}
\usepackage{pstcol}
\usepackage{indentfirst}
\usepackage{amsbsy}
\usepackage{times}
\usepackage[colorlinks=black]{hyperref}

\newtheorem{theorem}{Theorem}[section]
\newtheorem{defff}[theorem]{Definition}

\newtheorem{rem}[theorem]{Remark}
\newtheorem{lee}[theorem]{Lemma}

\numberwithin{equation}{section}

\begin{document}

\title{Existence and stability of global large strong solutions for the
Hall-MHD system \vspace{0.5cm}}
\author{{Maicon J. Benvenutti}{\thanks{%
Universidade Estadual de Campinas, IMECC-Departamento de Matem\'{a}tica, CEP
13083-859, Campinas-SP, Brazil. Email: mbenvenutti@hotmail.com. MJB was
supported by FAPESP, Brazil.}} \ \ \& \ {Lucas C. F. Ferreira}{\thanks{%
Universidade Estadual de Campinas, IMECC-Departamento de Matem\'{a}tica, CEP
13083-859, Campinas-SP, Brazil. Email: lcff@ime.unicamp.br. LCFF was
supported by FAPESP and CNPQ, Brazil. (corresponding author)}}}
\date{}
\maketitle
\date{}

\begin{abstract}
We consider the 3D incompressible Hall-MHD system and prove a stability
theorem for global large solutions under a suitable integrable hypothesis in
which one of the parcels is linked to the Hall term. As a byproduct, a class
of global strong solutions is obtained with large velocities and small initial
magnetic fields. Moreover, we prove the local-in-time well-posedness of $%
H^{2} $-strong solutions which improves previous regularity conditions on
initial data.

\

\noindent\textbf{AMS 2010 MSC:} 35Q35, 76D03, 35B35, 76E25, 76W05 \vspace{%
0.2cm}

 \vspace{0.2cm}
\end{abstract}

%\author{{Maicon J. Benvenutti}{\thanks{M. Benvenutti was supported by FAPESP, Brazil.}}\\{\small Universidade Estadual de Campinas}\\{\small IMECC-Departamento de Matem\'{a}tica} \\{\small {\ CEP 13083-859, Campinas-SP, Brazil}}\\{\small \texttt{E-mail:mbenvenutti@hotmail.com}}\vspace{0.5cm}\\{{Lucas C. F. Ferreira} {\thanks{L. Ferreira was supported by FAPESP and CNPQ, Brazil. (corresponding author)}}}\\{\small Universidade Estadual de Campinas, Departamento de Matem\'{a}tica}\\{\small {\ CEP 13083-859, Campinas-SP, Brazil}}\\{\small \texttt{E-mail:lcff@ime.unicamp.br}}}

\pagestyle{myheadings} %\thispagestyle{plain}
\markboth{Existence and stability for the Hall-MHD
system\hspace{4.5cm}}{\hspace{5.5cm} M.J. Benvenutti \ \& \ L.C.F. Ferreira}

\section{\protect\bigskip Introduction}

\hspace{0.6cm}This paper is concerned with the 3D incompressible Hall-MHD
system
\begin{equation}
\left\{
\begin{array}{rclll}
\displaystyle\partial _{t}u+[u.\nabla ]u+\nabla p-(\nabla \times b)\times b
& = & \mu \Delta u & \mbox{in} & \,\,(x,t)\in \mathbb{R}^{3}\times \lbrack
0,\infty ); \\
\displaystyle\partial _{t}b-\nabla \times (u\times b)+\nabla \times ((\nabla
\times b)\times b) & = & \gamma \Delta b & \mbox{in} & \,\,(x,t)\in \mathbb{R%
}^{3}\times \lbrack 0,\infty ); \\
div\,u & = & 0 & \mbox{in} & \,\,(x,t)\in \mathbb{R}^{3}\times \lbrack
0,\infty ), \\
&  &  &  &
\end{array}%
\right.  \label{1}
\end{equation}%
where $u=(u_{1}(x,t),u_{2}(x,t),u_{3}(x,t))$ is the velocity field, $%
p=p(x,t) $ is the scalar pressure field, $%
b=(b_{1}(x,t),b_{2}(x,t),b_{3}(x,t))$ is the magnetic field induced by the
charged fluid, $\mu >0$ and $\gamma >0$ are respectively the viscosity and
resistivity coefficients, $[u.\nabla ]=\sum_{i=1}^{3}u_{i}\partial _{x_{i}}$
and the symbol $\times $ stands for the usual three-dimensional
cross-product. The density of the fluid is assumed to be one by
normalization.

The system (\ref{1}) has been studied in the physics literature for decades
(see e.g. \cite{Liu}, \cite{Lighthill} and their references) and has
application in a number of physical fields such as geo-dynamo \cite{Mininni}%
, neutron stars \cite{Shalybkov} and magnetic reconnection in plasmas \cite%
{Forbes}. The reader is referred to \cite{Liu} (see also \cite{Bit}) for a
deduction of (\ref{1}) from two-fluids model, as well as from kinetic model,
considering a generalized Ohm law. In comparison with the usual
incompressible MHD system (see \cite{Sermange}), we have the new term $%
\nabla \times ((\nabla \times b)\times b)$ which is due to Hall effect and
prevents straightforward adaptations from arguments used in the mathematical
analysis of Navier-Stokes and related models.

Unlike MHD system that has an extensive variety of studies in classical
subjects such as existence of solutions, regularity criteria and stability
(see e.g. \cite{Miao-1}, \cite{Duvaut}, \cite{Fer-Vil-1}, \cite{Sermange},
\cite{Wu-1}, \cite{Zhou-1} and references therein), the influence of the
Hall term has been little explored on these topics. Indeed, Hall-MHD has
appeared only recently in the mathematical literature and there are
relatively a few works with this type of approach which are reviewed in what
follows. In \cite{Liu}, by using Galerkin's method, global in time existence
of weak solutions is proved in the periodic setting $L^{2}([0,1]^{3})$ for
the resistive ($\gamma >0$) and viscous case ($\mu >0$). The uniqueness of
weak solutions is still an open problem. Considering $\mu \geq 0$ and $%
\gamma >0,$ the authors of \cite{DC1} obtained, via energy method,
local-in-time well-posedness of strong solutions in $H^{m}(\mathbb{R}^{3})$
with $m>\frac{5}{2}$ as well as global well-posedness under small
conditions. They also showed blow-up criteria of first type for strong
solutions and a Liouville theorem for smooth stationary solutions. The main
point in \cite{DC1} was the control of the Hall term via diffusion induced
by the resistivity (see more details in the next paragraph). In \cite{DC2}
some blow-up criteria are studied and it is obtained a global well-posedness
result for small initial data in terms of Besov norm which can be considered
optimal in a suitable way that takes into account the scaling property for
the system with null velocity. A subclass of global strong axisymmetric
solutions was obtained in \cite{Fan}. By employing Fourier splitting method,
time-decay of Sobolev norms is showed in \cite{DC3} for a class of weak
solutions. A version of (\ref{1}) with magnetic fractional diffusion $%
(-\Delta )^{\alpha }$ was considered in \cite{DC4}, where it was proved
local well-posedness in Sobolev spaces for any $\alpha >\frac{2}{3}$ by
using the smoothing effects of the dissipation and local bounds for the
Sobolev norms through a multi-stage process. Regularity criteria for the
density-dependent case is studied in \cite{Fan2}. In \cite{DC5} it is shown
that the non-resistive system ($\gamma =0$) is not globally well-posed in
any Sobolev space $H^{m}(\mathbb{R}^{3})$ with $m>\frac{7}{2}$ in the sense
that either it is locally ill-posed or it is locally well-posed but there
exists an axisymmetric solution that loses the initial regularity in finite
time.

Due to the spatial derivative of high-order in a nonlinear term, the
Hall-MHD leads us to deal with higher regularity in the energy estimates
(see e.g. (\ref{aux-Hall-1})-(\ref{aux-Hall-3}) and (\ref{aux-Hall-4})-(\ref%
{aux-Hall-10}) in Section \ref{AAA2}), which introduces further difficulties
in handling the system. For comparison, let us recall briefly about the
issue of well-posedness for the incompressible Navier-Stokes and Euler
equations: the local well-posedness with large data and the global one with
small initial data in $H^{1}(\mathbb{R}^{3})$ for the Navier-Stokes
equations are obtained via an energy inequality where the nonlinearity, which
is of first order ($[u.\nabla ]u$), is estimated by using
Gagliardo-Nirenberg inequality and the diffusion controls the generated
second-order derivative (see \cite{Temam}). In the inviscid case, the local
well-posedness is obtained only in $H^{m}(\mathbb{R}^{3})$ with $m>\frac{5}{2%
}$ where the key inclusion $H^{m}(\mathbb{R}^{3})\hookrightarrow W^{1,\infty
}(\mathbb{R}^{3}) $ holds true (see \cite{Bertozzi}). In \cite{DC1}, Chae
\textit{et al.} mixed these two approaches to prove a local existence
theorem in $H^{m}(\mathbb{R}^{3})$, for $\gamma >0$ and $m>\frac{5}{2}$.
Another structural difference is that the second-order derivatives in the
Hall term seem to obstruct the parabolic regularization effect exhibited for
instance by MHD and Navier-Stokes systems.

Our first result improves (in the viscosity case) the one of \cite{DC1} by
proving local-in-time well-posedness of (\ref{1}) in $H^{2}(\mathbb{R}^{3}),$
and global well-posedness for small $H^{2}$-initial data (see Theorem \ref%
{a1}). Here we use accurate energy estimates and also the particular
structure of the Hall term.

Global existence of strong solutions of (\ref{1}) for large initial data is
still an open challenging problem. With respect to this matter, as far as we
know, there are just the above mentioned class of $2\frac{1}{2}$ dimensional
solutions of the form $(b,u)=(b(r,z)e^{\theta },u(r,z)e^{r}+u(r,z)e^{z})$ as
proved in \cite{Fan}. We observe that the two-dimensional symmetry is not
tractable due to the fact that in this situation the Hall term has just the
third component nonzero. Despite the helical symmetry is conserved for the
system, it is an open question to prove that they are global in time. Let us
again make a comparison with the Navier-Stokes and other classical systems.
For Navier-Stokes equations, there are global strong solutions in the
two-dimensional case (see e.g. \cite{Temam}), under the condition of axial
symmetry without swirl \cite{Ladyzhenskaya}, and in the presence of helical
symmetry \cite{Mahalov}. In \cite{T}, Ponce \textit{et al.} proved that
global solutions with a suitable property are stable in the sense that
solutions close to them are global as well. Fortunately, symmetric and
two-dimensional solutions satisfy the hypothesis required and this gives a
class of global large solutions which are genuinely three-dimensional
(although approximately symmetric or two-dimensional). Related results can
be found in \cite{Auscher}, \cite{Bardos}, \cite{Chemin}, \cite{Gallagher},
\cite{Gui}, \cite{Iftimie}, \cite{Mucha} and \cite{Rusin}. There are similar
theorems for inhomogeneous Navier-Stokes equations \cite{Abidi, Cai},
Boussinesq system \cite{Li2, Liu2} and MHD system \cite{Li}.

In this paper we extend the stability result of \cite{T} to the system (\ref%
{1}) (see Theorem \ref{theo2}). Again, the main difficult is the Hall term
that requires estimates to deal with higher derivatives in the nonlinear
term. Considering the global solutions obtained in \cite{T} for
Navier-Stokes equations, our stability result provides a class of global
strong solutions $(v,h)$ for (\ref{1}) with large velocities and small
initial magnetic fields (see Remark \ref{Rem33}).

This paper is organized in the following way: in Section \ref{AAA1} we give
some definitions, recall some basic inequalities and vector identities, and
discuss the formulation of the problem as well as the notions of weak and
strong solutions. Section \ref{AAA1-2} is devoted to state our results. In
Section \ref{AAA2} we obtain key estimates to deal with the system. Finally,
the results are proved in Section \ref{AAA3}.

\section{\protect\bigskip Preliminary}

\label{AAA1}

\subsection{\protect\bigskip Basic definitions and inequalities}

\hspace{0.6cm}Let us start with some basic definitions in order to formulate
the problem. We denote by $H^{m}(\mathbb{R}^{3})$ and $W^{m,p}(\mathbb{R}%
^{3})$ the usual three-dimensional vector Sobolev spaces (see \cite{Temam}).
The subscript $\sigma $ in $H^{m}(\mathbb{R}^{3})$ or in $L^{p}$ means that
the vector fields are divergence-free. The classical Helmholtz orthogonal
projection $\mathbb{P}$ onto the space of the solenoidal functions is
denoted by%
\begin{equation*}
\mathbb{P}:L^{2}(\mathbb{R}^{3})\longmapsto L_{\sigma }^{2}(\mathbb{R}^{3}).
\end{equation*}%
We recall the following particular cases of the Gagliardo-Nirenberg
inequality in $\mathbb{R}^{3}$ (see \cite{N})
\begin{equation}
\left\{
\begin{array}{rcll}
\Vert f\Vert _{L^{6}(\mathbb{R}^{3})} & \leq & C\Vert \nabla f\Vert _{L^{2}(%
\mathbb{R}^{3})}, & \forall f\in H^{1}(\mathbb{R}^{3}), \\
\Vert f\Vert _{L^{3}(\mathbb{R}^{3})} & \leq & C\Vert f\Vert _{L^{2}(\mathbb{%
R}^{3})}^{\frac{1}{2}}\Vert \nabla f\Vert _{L^{2}(\mathbb{R}^{3})}^{\frac{1}{%
2}}, & \forall f\in H^{1}(\mathbb{R}^{3}), \\
\Vert f\Vert _{L^{\infty }(\mathbb{R}^{3})} & \leq & C\Vert \nabla f\Vert
_{L^{2}(\mathbb{R}^{3})}^{\frac{1}{2}}\Vert \nabla ^{2}f\Vert _{L^{2}(%
\mathbb{R}^{3})}^{\frac{1}{2}}, & \forall f\in H^{2}(\mathbb{R}^{3}).%
\end{array}%
\right.  \label{gn}
\end{equation}%
Also, we will use some equivalent seminorms and norms that can be obtained
easily by Fourier transform. We have that
\begin{equation}
\left\{
\begin{array}{rclll}
\Vert \nabla ^{2}f\Vert _{L^{2}(\mathbb{R}^{3})}^{2} & \cong & \Vert \Delta
f\Vert _{L^{2}(\mathbb{R}^{3})}^{2}=\Vert \nabla \times \nabla \times f\Vert
_{L^{2}(\mathbb{R}^{3})}^{2}+\Vert \nabla div\,f\Vert _{L^{2}(\mathbb{R}%
^{3})}^{2} & \mbox{in} & H^{2}(\mathbb{R}^{3}), \\
\Vert \nabla ^{3}f\Vert _{L^{2}(\mathbb{R}^{3})}^{2} & \cong & \Vert
div\,\Delta f\Vert _{L^{2}(\mathbb{R}^{3})}^{2}+\Vert \nabla \times \Delta
f\Vert _{L^{2}(\mathbb{R}^{3})}^{2} & \mbox{in} & H^{3}(\mathbb{R}^{3}), \\
\Vert \nabla ^{3}f\Vert _{L^{2}(\mathbb{R}^{3})}^{2} & \cong & \Vert \nabla
\times \Delta f\Vert _{L^{2}(\mathbb{R}^{3})}^{2} & \mbox{in} & H_{\sigma
}^{3}(\mathbb{R}^{3}), \\
\Vert \nabla ^{3}f\Vert _{L^{2}(\mathbb{R}^{3})}^{2} & \cong & \Vert
div\,\Delta f\Vert _{L^{2}(\mathbb{R}^{3})}^{2} & \mbox{in} & H_{p}^{3}(%
\mathbb{R}^{3}),%
\end{array}%
\right.  \label{aux-equiv-1}
\end{equation}%
and
\begin{equation}
\left\{
\begin{array}{rcl}
\Vert f\Vert _{H^{2}(\mathbb{R}^{3})}^{2} & \cong & \Vert f\Vert _{L^{2}(%
\mathbb{R}^{3})}^{2}+\Vert \nabla f\Vert _{L^{2}(\mathbb{R}^{3})}^{2}+\Vert
\Delta f\Vert _{L^{2}(\mathbb{R}^{3})}^{2}, \\
\Vert f\Vert _{H^{3}(\mathbb{R}^{3})}^{2} & \cong & \Vert f\Vert _{L^{2}(%
\mathbb{R}^{3})}^{2}+\Vert \nabla f\Vert _{L^{2}(\mathbb{R}^{3})}^{2}+\Vert
\Delta f\Vert _{L^{2}(\mathbb{R}^{3})}^{2}+\Vert div\,\Delta f\Vert _{L^{2}(%
\mathbb{R}^{3})}^{2}+\Vert \nabla \times \Delta f\Vert _{L^{2}(\mathbb{R}%
^{3})}^{2}, \\
\Vert f\Vert _{H_{\sigma }^{3}(\mathbb{R}^{3})}^{2} & \cong & \Vert f\Vert
_{L^{2}(\mathbb{R}^{3})}^{2}+\Vert \nabla f\Vert _{L^{2}(\mathbb{R}%
^{3})}^{2}+\Vert \Delta f\Vert _{L^{2}(\mathbb{R}^{3})}^{2}+\Vert \nabla
\times \Delta f\Vert _{L^{2}(\mathbb{R}^{3})}^{2}, \\
\Vert f\Vert _{H_{p}^{3}(\mathbb{R}^{3})}^{2} & \cong & \Vert f\Vert _{L^{2}(%
\mathbb{R}^{3})}^{2}+\Vert \nabla f\Vert _{L^{2}(\mathbb{R}^{3})}^{2}+\Vert
\Delta f\Vert _{L^{2}(\mathbb{R}^{3})}^{2}+\Vert div\,\Delta f\Vert _{L^{2}(%
\mathbb{R}^{3})}^{2}.%
\end{array}%
\right.  \label{aux-equiv-2}
\end{equation}

\subsection{\protect\bigskip Vector identities}

\hspace{0.6cm}Here we recall some vector equalities which will be useful in
order to deal with the Hall term. We have (see \cite{M})
\begin{equation*}
\left\{
\begin{array}{rcl}
\Delta A & = & \nabla div\,A-\nabla \times \nabla \times A, \\
(\nabla \times A)\times B & = & [A.\nabla ]B+[B.\nabla ]A+A\times (\nabla
\times B)-\nabla (A.B), \\
\nabla \times (A\times B) & = & A(div\,B)-B(div\,A)+[B.\nabla ]A-[A.\nabla
]B,%
\end{array}%
\right.
\end{equation*}%
from which we obtain

\begin{eqnarray}
\nabla \times ((\nabla \times A)\times B)-(\nabla \times \nabla \times
A)\times B &=&(\nabla \times A)(div\,B)-2[(\nabla \times A).\nabla ]B
\label{345-1} \\
&&-(\nabla \times A)\times (\nabla \times B)+\nabla ((\nabla \times A).B)
\notag
\end{eqnarray}%
and

\begin{eqnarray}
\nabla \times ((\nabla \times \nabla \times A)\times B)-(\nabla \times
\nabla \times \nabla \times A)\times B &=&(\nabla \times \nabla \times
A)(div\,B)  \notag \\
&&-2[(\nabla \times \nabla \times A).\nabla ]B  \notag \\
&&-(\nabla \times \nabla \times A)\times (\nabla \times B)  \label{345-2} \\
&&+\nabla ((\nabla \times \nabla \times A).B).  \notag
\end{eqnarray}

\subsection{\protect\bigskip Weak and strong solutions}

\hspace{0.6cm}Consider the operators
\begin{align}
& A_{1}:H_{\sigma }^{1}(\mathbb{R}^{3})\longmapsto H_{\sigma }^{-1}(\mathbb{R%
}^{3})\mbox{ defined by }\left\langle A_{1}[u],v\right\rangle =\mu \int_{%
\mathbb{R}^{2}}\nabla u.\nabla v\,dx;  \notag \\
& A_{2}:H_{\sigma }^{1}(\mathbb{R}^{3})\times H_{\sigma }^{1}(\mathbb{R}%
^{3})\longmapsto H_{\sigma }^{-1}(\mathbb{R}^{3})\mbox{ defined by }%
\left\langle A_{2}[u,h],v\right\rangle =\int_{\mathbb{R}^{2}}([u.\nabla
]h).v\,dx;  \notag \\
& A_{3}:H_{{}}^{1}(\mathbb{R}^{3})\times H_{{}}^{1}(\mathbb{R}%
^{3})\longmapsto H_{\sigma }^{-1}(\mathbb{R}^{3})\mbox{ defined by }%
\left\langle A_{3}[b,h],v\right\rangle =-\int_{\mathbb{R}^{2}}((\nabla
\times b)\times h).v\,dx;  \notag \\
& B_{1}:H^{1}(\mathbb{R}^{3})\longmapsto H^{-2}(\mathbb{R}^{3})%
\mbox{
defined by }\left\langle B_{1}[b],w\right\rangle =\gamma \int_{\mathbb{R}%
^{2}}\nabla b.\nabla w\,dx;  \notag \\
& B_{2}:H_{\sigma }^{1}(\mathbb{R}^{3})\times H^{1}(\mathbb{R}%
^{3})\longmapsto H^{-2}(\mathbb{R}^{3})\mbox{ defined by }\left\langle
B_{2}[u,b],w\right\rangle =-\int_{\mathbb{R}^{2}}(\nabla \times (u\times
b)).w\,dx;  \notag \\
& B_{3}:H^{1}(\mathbb{R}^{3})\times H^{1}(\mathbb{R}^{3})\longmapsto H^{-2}(%
\mathbb{R}^{3})\mbox{ defined by }\left\langle B_{3}[b,h],w\right\rangle
=\int_{\mathbb{R}^{2}}((\nabla \times b)\times h).(\nabla \times w)\,dx.
\notag
\end{align}

It is straightforward to prove by Gagliardo-Nirenberg type inequalities that
these operators are well-defined and continuous. We consider the usual weak
formulation for (\ref{1})
\begin{equation}
\left\{
\begin{array}{rcl}
\frac{d}{dt}u+A_{1}[u]+A_{2}[u,u]+A_{3}[b,b] & = & 0\mbox{ in }%
L^{1}((0,T),H_{\sigma }^{-1}(\mathbb{R}^{3})); \\
\frac{d}{dt}b+B_{1}[b]+B_{2}[u,b]+B_{3}[b,b] & = & 0\mbox{ in }%
L^{1}((0,T),H^{-2}(\mathbb{R}^{3})); \\
\mathbb{P}[u] & = & u.%
\end{array}%
\right.  \label{3}
\end{equation}%
For $(u_{0},b_{0})\in L_{\sigma }^{2}(\mathbb{R}^{3})\times L^{2}(\mathbb{R}%
^{3})$ and $0<T<\infty ,$ we say that $(u,b)$ is a weak solution in $(0,T)$
for (\ref{1}) with initial data $(u_{0},b_{0})$ if
\begin{equation}
\left\{
\begin{array}{rcl}
u & \in & L^{2}((0,T);H_{\sigma }^{1}(\mathbb{R}^{3}))\cap L^{\infty
}((0,T);L_{\sigma }^{2}(\mathbb{R}^{3})), \\
b & \in & L^{2}((0,T);H^{1}(\mathbb{R}^{3}))\cap L^{\infty }((0,T);L^{2}(%
\mathbb{R}^{3})), \\
\frac{d}{dt}u & \in & L^{1}((0,T),H_{\sigma }^{-1}(\mathbb{R}^{3})), \\
\frac{d}{dt}b & \in & L^{1}((0,T),H^{-2}(\mathbb{R}^{3}))%
\end{array}%
\right.  \label{aux-weak-1}
\end{equation}%
and $(u,b)$ satisfies (\ref{3}). In the case $(0,\infty )$ (global
solutions), we assume that $(u,b)$ satisfies (\ref{3}) and (\ref{aux-weak-1}%
) for all $0<T<\infty .$

\begin{rem}
\label{Rem1}If $(u,b)$ is a weak solution, then $(u,b)\in C_{w}([0,T],L^{2}(%
\mathbb{R}^{3}))$ and the initial data condition is satisfied in an
appropriate sense of weak limit (see \cite{Temam}).
\end{rem}

Inspired on the classical mathematical literature, it is natural to consider
class of solutions in spaces where energy estimates provide, at least, local
well-posedness. For Navier-Stokes equations (and also MHD), the space $%
L^{2}((0,T);H^{2})\cap L^{\infty }((0,T);H^{1})$ is commonly used together
with $H^{1}$ initial data. These solutions are strong in the sense that they
have $H^{1}$-continuous orbits (i.e., belong to $C([0,T),H^{1})$) and
satisfy their respective systems in $L^{2}$ for almost everywhere $t\in
(0,T) $.

Due the second-order derivative in the non-linear part of (\ref{1}), the
above space is not appropriated to perform suitable energy estimates.
However, using the special structure of Hall term $\nabla \times ((\nabla
\times b)\times b)$, we will prove the local well-posedness in $%
L^{2}((0,T);H_{\sigma }^{3}\times H^{3})\cap L^{\infty }((0,T);H_{\sigma
}^{2}\times H^{2})$ for $H^{2}$ initial data and these solutions have $H^{2}$%
-continuous orbits and satisfy the system in $L^{2}$, for almost everywhere $%
t\in (0,T)$. So, we establish the following definition.

\begin{defff}
\label{Def-1}(Strong solution) Let $(u_{0},b_{0})\in H_{\sigma }^{2}(\mathbb{%
R}^{3})\times H^{2}(\mathbb{R}^{3}).$ For $0<T<\infty ,$ we say that $(u,b)$
is a strong solution in $(0,T)$ for (\ref{1}) with initial data $%
(u_{0},b_{0})$ if $(u,b)$ verifies (\ref{3}) and belongs to class
\begin{equation}
\left\{
\begin{array}{rcl}
u & \in & L^{2}((0,T);H_{\sigma }^{3}(\mathbb{R}^{3}))\cap L^{\infty
}((0,T);H_{\sigma }^{2}(\mathbb{R}^{3})), \\
b & \in & L^{2}((0,T);H^{3}(\mathbb{R}^{3}))\cap L^{\infty }((0,T);H^{2}(%
\mathbb{R}^{3})).%
\end{array}%
\right.  \label{aux-reg-1}
\end{equation}%
In the case of $(0,\infty )$ (global solutions), we assume that $(u,b)$
satisfies (\ref{3}) and (\ref{aux-reg-1}) for all $0<T<\infty .$
\end{defff}

\begin{rem}
Indeed, we are going to prove uniqueness of weak solutions in a class larger
than (\ref{aux-reg-1}), namely $L^{4}((0,T),H_{\sigma }^{1}(\mathbb{R}%
^{3}))\times L^{4}((0,T),H^{2}(\mathbb{R}^{3}))$ (see Theorem \ref{a1}).
\end{rem}

\begin{rem}
\label{Rem3} It is straightforward to check that a strong solution of (\ref%
{1}) satisfies
\begin{equation}
\frac{d}{dt}\Delta u\in L^{2}((0,T);H_{\sigma }^{-1}(\mathbb{R}^{3})) \ \
\text{and} \ \ \frac{d}{dt}\Delta b\in L^{2}((0,T);H^{-1}(\mathbb{R}^{3})).
\end{equation}
Therefore
\begin{equation}
u,b\in C([0,T),H^{2}(\mathbb{R}^{3})),  \label{aux-continuity}
\end{equation}%
$\left\langle \frac{d}{dt}\Delta u,\Delta u\right\rangle =\frac{1}{2}\frac{d%
}{dt}\Vert \Delta u\Vert _{L^{2}(\mathbb{R}^{3})}^{2}$ and $\left\langle
\frac{d}{dt}\Delta b,\Delta b\right\rangle =\frac{1}{2}\frac{d}{dt}\Vert
\Delta b\Vert _{L^{2}(\mathbb{R}^{3})}^{2}$ (see \cite{Temam} for further
details). The same is obviously true for spatial derivatives of lower order.
\newline
%%%%%%%%%%%%%%%%%%%%%%%%%%%%%%%%%%%%%%%%%%%%%%%%%%%%%%%%%%%%%%%%%%%%%%%%%%%%%%%%%%%%%%%%%%%%%%%%%%%%%%%%%
%%%%%%%%%%%%%%%%%%%%%%%%%%%%%%%%%%%%%%%%%%%%%%%%%%%%%%%%%%%%%%%%%%%%%%%%%%%%%%%%%%%%%%%%%%%%%%%%%%%%%%%%%
\end{rem}

\section{\protect\bigskip Results}

\label{AAA1-2}

\hspace{0.6cm}In this section we state our results. We start with a result
which improves the initial data regularity condition in \cite{DC1} for local-in-time well-posedness.

\subsection{\protect\bigskip Local-in-time well-posedness in $H^{2}$}

\begin{theorem}
\label{a1} Let $(u_{0},b_{0})\in H_{\sigma }^{2}(\mathbb{R}^{3})\times H^{2}(%
\mathbb{R}^{3})$. Then, there exist $T=T(\Vert u_{0}\Vert _{H_{\sigma }^{2}(%
\mathbb{R}^{3})},\Vert B_{0}\Vert _{H^{2}(\mathbb{R}^{3})})>0$ and a strong
solution $(u,b)$ of (\ref{1}) in $(0,T)$ with initial data $(u_{0},b_{0}).$
This solution is the unique weak solution in $L^{4}((0,T),H_{\sigma }^{1}(%
\mathbb{R}^{3}))\times L^{4}((0,T),H^{2}(\mathbb{R}^{3}))$. Furthermore, if $%
\Vert u_{0}\Vert _{H_{\sigma }^{2}(\mathbb{R}^{3})}^{2}+\Vert b_{0}\Vert
_{H^{2}(\mathbb{R}^{3})}^{2}$ is small enough, then the solution is global
in time. Finally, if $T<\infty $ is the maximal existence time, then
\begin{equation}
\int_{0}^{T}\left( \Vert \nabla u\Vert _{L^{2}(\mathbb{R}^{3})}^{4}+\Vert
\nabla b\Vert _{L^{2}(\mathbb{R}^{3})}^{4}+\Vert \Delta b\Vert _{L^{2}(%
\mathbb{R}^{3})}^{4}\right) dt=\infty .  \label{asd}
\end{equation}
\end{theorem}

\subsection{\protect\bigskip Global stability of large solutions}

\hspace{0.6cm}In the next theorem we obtain stability of large global strong
solutions whose integral in (\ref{asd}) is finite with $T=\infty $. Notice that this condition is natural because we are dealing with global solutions.

\begin{theorem}
\label{theo2} Let $(u,b)$ be a global strong solution of (\ref{1}) with
initial data $(u_{0},b_{0})\in H_{\sigma }^{2}(\mathbb{R}^{3})\times H^{2}(%
\mathbb{R}^{3})$ and satisfying
\begin{equation}
\int_{0}^{\infty }\left( \Vert \nabla u\Vert _{L^{2}(\mathbb{R}%
^{3})}^{4}+\Vert \nabla b\Vert _{L^{2}(\mathbb{R}^{3})}^{4}+\Vert \Delta
b\Vert _{L^{2}(\mathbb{R}^{3})}^{4}\right) dt<\infty .  \label{444}
\end{equation}%
There exists $\delta >0$ such that if $(v_{0},h_{0})\in H_{\sigma }^{2}(%
\mathbb{R}^{3})\times H^{2}(\mathbb{R}^{3})$ and
\begin{equation}
\Vert u_{0}-v_{0}\Vert _{H_{\sigma }^{2}(\mathbb{R}^{3})}^{2}+\Vert
b_{0}-h_{0}\Vert _{H^{2}(\mathbb{R}^{3})}^{2}<\delta ,  \label{aux-stab-1}
\end{equation}%
then the strong solution $(v,h)$ with initial data $(v_{0},h_{0})$ is global
in time. Furthermore, there exists $M=M(\delta )$ with $M(\delta )\overset{%
\delta \rightarrow 0}{\longrightarrow }0$ such that
\begin{equation*}
\sup_{t\geq 0}\left( \Vert u(t)-v(t)\Vert _{H_{\sigma }^{2}(\mathbb{R}%
^{3})}^{2}+\Vert b(t)-h(t)\Vert _{H^{2}(\mathbb{R}^{3})}^{2}\right) \leq
M(\delta ).
\end{equation*}
\end{theorem}

\begin{rem}
\label{Rem33} In Theorem \ref{a1}, global strong solutions are obtained for
small initial velocities and magnetic fields. We can use Theorem \ref{theo2}
to provide a class of global strong solutions with large initial velocities
and small initial magnetic fields. Let us consider the classical
incompressible Navier-Stokes equations
\begin{equation}
\left\{
\begin{array}{rclll}
\displaystyle\partial _{t}u+[u.\nabla ]u+\nabla p & = & \mu \Delta u & %
\mbox{in} & \,\,(x,t)\in \mathbb{R}^{3}\times \lbrack 0,\infty ); \\
div\,u & = & 0 & \mbox{in} & \,\,(x,t)\in \mathbb{R}^{3}\times \lbrack
0,\infty ).%
\end{array}%
\right.   \label{1112}
\end{equation}%
As pointed out in Introduction, the paper \cite{T} provides a class of
global large solutions for (\ref{1112}) satisfying
\begin{equation}
\int_{0}^{\infty }\Vert \nabla u(s)\Vert _{L^{2}(\mathbb{R}%
^{3})}^{4}ds<\infty .  \label{cond1}
\end{equation}%
We have that $(u,b)\equiv (u,0)$ is a global strong solution for (\ref{1})
and verifies (\ref{444}). If $(v,h)$ is a local-in-time strong solution for (\ref{1}) (given by Theorem \ref{a1}) such that $v(0)$ is close to $u(0)$ and $h(0)$ is small enough, then $(v,h)$ is
also a global strong solution.
\end{rem}

%%%%%%%%%%%%%%%%%%%%%%%%%%%%%%%%%%%%%%%%%%%%%%%%%%%%%%%%%%%%%%%%%%%%%%%%%%%%%%%%%%%%%%%%%%%%%%%%%%%%%%%%%
%%%%%%%%%%%%%%%%%%%%%%%%%%%%%%%%%%%%%%%%%%%%%%%%%%%%%%%%%%%%%%%%%%%%%%%%%%%%%%%%%%%%%%%%%%%%%%%%%%%%%%%%%

\section{\protect\bigskip Key estimates}

\label{AAA2}

\hspace{0.6cm}We start with two lemmas which contain energy estimates that
will be used to prove the results stated in Section \ref{AAA1-2}. For the
sake of presentation, the proof of Lemma \ref{10} is postponed for
Subsection \ref{proof-lem10}.

\begin{lee}
\label{10} Let $(u,b)$ be a strong solution of (\ref{1}) in $(0,T)$
according to Definition \ref{Def-1}. We have that
\begin{equation}
\frac{1}{2}\frac{d}{dt}\left( \Vert u(t)\Vert _{L^{2}}^{2}+\Vert b(t)\Vert
_{L^{2}}^{2}\right) +\mu \Vert \nabla u(t)\Vert _{L^{2}}^{2}+\gamma \Vert
\nabla b(t)\Vert _{L^{2}}^{2}=0,\,\,\forall \,0\leq t<T.  \label{in1}
\end{equation}%
Furthermore, there are constants $C_{0}=C_{0}(\mu ,\,\gamma )>0$ and $%
C_{1}=C_{1}(\mu ,\,\gamma )>0$ such that
\begin{align}
\frac{1}{2}\frac{d}{dt}\left( \Vert \nabla u(t)\Vert _{L^{2}}^{2}\right. &
+\left. \Vert \nabla b(t)\Vert _{L^{2}}^{2}\right) +\frac{\mu }{2}\Vert
\Delta u(t)\Vert _{L^{2}}^{2}+\frac{\gamma }{2}\Vert \Delta b(t)\Vert
_{L^{2}}^{2}  \notag \\
& \leq C_{0}\left( \Vert \nabla u(t)\Vert _{L^{2}}^{2}+\Vert \nabla
b(t)\Vert _{L^{2}}^{2}\right) ^{3}+\Vert \Delta b(t)\Vert _{L^{2}}^{4}\Vert
\nabla b(t)\Vert _{L^{2}}^{2},\,\,\forall \,0\leq t<T,  \label{in2}
\end{align}%
and
\begin{align}
\frac{1}{2}\frac{d}{dt}\left( \Vert \Delta u(t)\Vert _{L^{2}}^{2}\right. &
+\left. \Vert \Delta b(t)\Vert _{L^{2}}^{2}\right) +\frac{\mu }{4}\Vert
\nabla \times \Delta u(t)\Vert _{L^{2}}^{2}+\frac{\gamma }{4}\left( \Vert
\nabla \times \Delta b(t)\Vert _{L^{2}}^{2}+\Vert div\,\Delta b(t)\Vert
_{L^{2}}^{2}\right)  \notag \\
& \leq \frac{\gamma }{4}\Vert \Delta b(t)\Vert _{L^{2}}^{2}+C_{1}\left(
\Vert \nabla u(t)\Vert _{L^{2}}^{2}+\Vert \nabla b(t)\Vert
_{L^{2}}^{2}+\Vert \Delta b(t)\Vert _{L^{2}}^{2}\right) ^{3}  \notag \\
& +C_{1}\Vert \Delta u(t)\Vert _{L^{2}}^{2}\left( \Vert \nabla u(t)\Vert
_{L^{2}}^{4}+\Vert \nabla b(t)\Vert _{L^{2}}^{4}\right) ,\,\,\forall \,0\leq
t<T.  \label{in4}
\end{align}

\begin{rem}
\label{8}Using (\ref{in1})-(\ref{in4}) together with the equivalent norms (%
\ref{aux-equiv-1})-(\ref{aux-equiv-2}), we have that the above solution
satisfies
\begin{align}
& \frac{1}{2}\frac{d}{dt}\left( \Vert u(t)\Vert _{H^{2}(\mathbb{R}%
^{3})}^{2}+\Vert b(t)\Vert _{H^{2}(\mathbb{R}^{3})}^{2}\right) +\frac{\mu }{4%
}\left( \Vert \nabla u(t)\Vert _{L^{2}(\mathbb{R}^{3})}^{2}+\Vert \Delta
u(t)\Vert _{L^{2}(\mathbb{R}^{3})}^{2}+\Vert \nabla \times \Delta u(t)\Vert
_{L^{2}(\mathbb{R}^{3})}^{2}\right)  \notag \\
& +\frac{\gamma }{4}\left( \Vert \nabla b(t)\Vert _{L^{2}(\mathbb{R}%
^{3})}^{2}+\Vert \Delta b(t)\Vert _{L^{2}(\mathbb{R}^{3})}^{2}+\Vert \nabla
\times \Delta b(t)\Vert _{L^{2}(\mathbb{R}^{3})}^{2}+\Vert div\,\Delta
b(t)\Vert _{L^{2}(\mathbb{R}^{3})}^{2}\right)  \notag \\
& \leq C\left( \Vert u(t)\Vert _{H^{2}(\mathbb{R}^{3})}^{2}+\Vert b(t)\Vert
_{H^{2}(\mathbb{R}^{3})}^{2}\right) \left( \Vert \nabla u(t)\Vert _{L^{2}(%
\mathbb{R}^{3})}^{4}+\Vert \nabla b(t)\Vert _{L^{2}(\mathbb{R}%
^{3})}^{4}+\Vert \Delta b(t)\Vert _{L^{2}(\mathbb{R}^{3})}^{4}\right) .
\label{cach}
\end{align}
\end{rem}
\end{lee}

\bigskip

The subject of the next lemma is to show that solutions as in Lemma \ref{10}
under the condition (\ref{444}) satisfy a stronger estimate that gives some
control on the Hall-term.

\begin{lee}
\label{9} Let $(u,b)$ be a global strong solution that satisfies (\ref{444}%
). Then
\begin{equation}
\int_{0}^{\infty }\Vert \nabla u(s)\Vert _{L^{2}}^{4}+\Vert \nabla b(s)\Vert
_{L^{2}}^{4}+\Vert \Delta u(s)\Vert _{L^{2}}^{4}+\Vert \Delta b(s)\Vert
_{L^{2}}^{4}+\Vert \nabla \times \Delta b(s)\Vert _{L^{2}}^{2}+\Vert
div\,\Delta b(s)\Vert _{L^{2}}^{2}\,ds<\infty .  \label{aux-ineq-2}
\end{equation}
\end{lee}

\textbf{Proof. }Let $f(t)=\Vert u(t)\Vert _{H^{2}(\mathbb{R}^{3})}^{2}+\Vert
b(t)\Vert _{H^{2}(\mathbb{R}^{3})}^{2}$. So, by (\ref{cach}), we have
\begin{equation}
\frac{1}{2}\frac{d}{dt}f(t)\leq Cf(t)\left( \Vert \nabla u(s)\Vert _{L^{2}(%
\mathbb{R}^{3})}^{4}+\Vert \nabla b(s)\Vert _{L^{2}(\mathbb{R}%
^{3})}^{4}+\Vert \Delta b(s)\Vert _{L^{2}(\mathbb{R}^{3})}^{4}\right) .
\label{aux-ineq-1}
\end{equation}%
If we define
\begin{equation*}
M=\int_{0}^{\infty }\left( \Vert \nabla u(s)\Vert _{L^{2}(\mathbb{R}%
^{3})}^{4}+\Vert \nabla b(s)\Vert _{L^{2}(\mathbb{R}^{3})}^{4}+\Vert \Delta
b(s)\Vert _{L^{2}(\mathbb{R}^{3})}^{4}\right) ds<\infty
\end{equation*}%
and use Gronwall inequality in (\ref{aux-ineq-1}), we obtain
\begin{equation}
\Vert u(t)\Vert _{H^{2}(\mathbb{R}^{3})}^{2}+\Vert b(t)\Vert _{H^{2}(\mathbb{%
R}^{3})}^{2}\leq e^{2CM}\left( \Vert u(0)\Vert _{H^{2}(\mathbb{R}%
^{3})}^{2}+\Vert b(0)\Vert _{H^{2}(\mathbb{R}^{3})}^{2}\right) .
\label{ineq-3}
\end{equation}%
On the other side, integrating (\ref{cach}), we also get
\begin{align}
\frac{\gamma }{4}\int_{0}^{t}\left( \Vert \nabla \times \Delta b(s)\Vert
_{L^{2}(\mathbb{R}^{3})}^{2}\right. & +\left. \Vert div\,\Delta b(s)\Vert
_{L^{2}(\mathbb{R}^{3})}^{2}\right) ds+\frac{\mu }{4}\int_{0}^{t}\Vert
\Delta u(s)\Vert _{L^{2}(\mathbb{R}^{3})}^{2}ds  \notag \\
& \leq CM\sup_{s>0}\left\{ \Vert u(s)\Vert _{H^{2}(\mathbb{R}%
^{3})}^{2}+\Vert b(s)\Vert _{H^{2}(\mathbb{R}^{3})}^{2}\right\} .
\label{aux-ineq-4}
\end{align}%
The condition (\ref{444}) and inequalities (\ref{ineq-3})-(\ref{aux-ineq-4})
give (\ref{aux-ineq-2}).

\begin{flushright}
\rule{2mm}{2mm}
\end{flushright}

In order to obtain the stability result, we need to estimate the difference
between two strong solutions of (\ref{1}). For this, we have two lemmas
whose proofs are relatively long and so we also postpone them for later (see
subsections \ref{proof-lem11} and \ref{proof-lem11-b}).

\begin{lee}
\label{11} Let $(u,b)$ and $(v,h)$ be strong solutions of (\ref{1}). If $%
U=v-u$ and $B=h-b$ then
\begin{equation}
\left\{
\begin{array}{rcl}
\displaystyle\partial _{t}U-\mu \Delta U & = & -P[[U.\nabla ]U]-P[[U.\nabla
]u]-P[[u.\nabla ]U]+P[(\nabla \times B)\times B] \\
& + & P[(\nabla \times B)\times b]+P[(\nabla \times b)\times B]; \\
&  &  \\
\displaystyle\partial _{t}B-\gamma \Delta B & = & \nabla \times (U\times
B)+\nabla \times (U\times b)+\nabla \times (u\times B)-\nabla \times
((\nabla \times B)\times B) \\
& - & \nabla \times ((\nabla \times B)\times b)-\nabla \times ((\nabla
\times b)\times B); \\
&  &  \\
div\,U & = & 0.%
\end{array}%
\right.  \label{4}
\end{equation}%
Furthermore, there are constants $C_{2}=C_{2}(\mu ,\,\gamma )>0$ and $%
C_{3}=C_{3}(\mu ,\,\gamma )>0$ such that
\begin{align}
& \frac{1}{2}\frac{d}{dt}\left( \Vert U(t)\Vert _{L^{2}}^{2}+\Vert B(t)\Vert
_{L^{2}}^{2}\right) +\frac{\mu }{2}\Vert \nabla U(t)\Vert _{L^{2}}^{2}+\frac{%
\gamma }{2}\Vert \nabla B(t)\Vert _{L^{2}}^{2}  \notag \\
& \leq C_{2}\left( \Vert \nabla u(t)\Vert _{L^{2}}^{4}+\Vert \nabla
b(t)\Vert _{L^{2}}^{4}+\Vert \Delta b(t)\Vert _{L^{2}}^{4}\right) (\Vert
U(t)\Vert _{L^{2}}^{2}+\left. \Vert B(t)\Vert _{L^{2}}^{2}\right)
,\,\,\forall \,0\leq t<T,  \label{in5}
\end{align}%
and%
\begin{align}
& \frac{1}{2}\frac{d}{dt}\left( \Vert \nabla U(t)\Vert _{L^{2}}^{2}+\Vert
\nabla B(t)\Vert _{L^{2}}^{2}\right) +\frac{\mu }{2}\Vert \Delta U(t)\Vert
_{L^{2}}^{2}+\frac{\gamma }{2}\Vert \Delta B(t)\Vert _{L^{2}}^{2}  \notag \\
& \leq C_{3}\left( \Vert \nabla U(t)\Vert _{L^{2}}^{2}+\Vert \nabla
B(t)\Vert _{L^{2}}^{2}+\Vert \Delta B(t)\Vert _{L^{2}}^{2}\right) (\Vert
\nabla u(t)\Vert _{L^{2}}^{4}+\Vert \nabla b(t)\Vert _{L^{2}}^{4}+\Vert
\Delta b(t)\Vert _{L^{2}}^{4})  \notag \\
& +C_{3}\left( \Vert \nabla U(t)\Vert _{L^{2}}^{2}+\Vert \nabla B(t)\Vert
_{L^{2}}^{2}+\Vert \Delta B(t)\Vert _{L^{2}}^{2}\right) ^{3},\,\,\forall
\,0\leq t<T.  \label{in50}
\end{align}
\end{lee}

\begin{rem}
\label{Rem-4-10}The inequality (\ref{in5}) holds true for two weak solutions
$(u,b)$ and $(v,h)$ belonging to $L^{4}((0,T),H_{\sigma }^{1}(\mathbb{R}%
^{3}))\times L^{4}((0,T),H^{2}(\mathbb{R}^{3}))$.
\end{rem}

\begin{lee}
\label{11-b} Let $(u,b),$ $(v,h)$ and $(U,B)$ as in Lemma \ref{11}. There is
a constant $C_{4}=C_{4}(\mu ,\,\gamma )>0$ such that
\begin{align}
& \frac{1}{2}\frac{d}{dt}\left( \Vert \Delta U\Vert _{L^{2}}^{2}+\Vert
\Delta B\Vert _{L^{2}}^{2}\right) +\frac{\mu }{4}\Vert \nabla \times \Delta
U\Vert _{L^{2}}^{2}+\frac{\gamma }{4}\left( \Vert \nabla \times \Delta
B\Vert _{L^{2}}^{2}+\Vert div\,\Delta B\Vert _{L^{2}}^{2}\right)  \notag \\
& \leq C_{4}\left( \Vert \nabla U\Vert _{L^{2}}^{2}+\Vert \nabla B\Vert
_{L^{2}}^{2}+\Vert \Delta U\Vert _{L^{2}}^{2}+\Vert \Delta B\Vert
_{L^{2}}^{2}\right) (\Vert \nabla u\Vert _{L^{2}}^{4}+\Vert \nabla b\Vert
_{L^{4}}^{4}+\Vert \Delta u\Vert _{L^{2}}^{4}+\Vert \Delta b\Vert
_{L^{2}}^{4})  \notag \\
& +\frac{\mu }{4}\Vert \Delta U\Vert _{L^{2}}^{2}+\frac{\gamma }{4}\Vert
\Delta B\Vert _{L^{2}}^{2}+C_{4}\left( \Vert \nabla B\Vert
_{L^{2}}^{2}+\Vert \Delta B\Vert _{L^{2}}^{2}\right) \left( \Vert \nabla
\times \Delta b\Vert _{L^{2}}^{2}+\Vert div\,\Delta b\Vert
_{L^{2}}^{2}\right)  \notag \\
& +C_{4}\left( \Vert \nabla U\Vert _{L^{2}}^{2}+\Vert \nabla B\Vert
_{L^{2}}^{2}+\Vert \Delta U\Vert _{L^{2}}^{2}+\Vert \Delta B\Vert
_{L^{2}}^{2}\right) ^{3},\,\,\forall \,0\leq t<T.  \label{in300}
\end{align}

\begin{rem}
\label{12} Let $U=v-u$, $B=h-b$ and
\begin{equation*}
L(t)=\Vert \nabla u\Vert _{L^{2}}^{4}+\Vert \nabla b\Vert _{L^{2}}^{4}+\Vert
\Delta u\Vert _{L^{2}}^{4}+\Vert \Delta b\Vert _{L^{2}}^{4}+\Vert \nabla
\times \Delta b\Vert _{L^{2}}^{2}+\Vert div\,\Delta b\Vert _{L^{2}}^{2}.
\end{equation*}%
It follows from (\ref{in5})-(\ref{in300}) that
\begin{align*}
\frac{d}{dt}\left( \Vert U\Vert _{H^{2}}^{2}\right. & +\left. \Vert B\Vert
_{H^{2}}^{2}\right) +\frac{\mu }{4}\left( \Vert \nabla U\Vert _{L^{2}(%
\mathbb{R}^{3})}^{2}+\Vert \Delta U\Vert _{L^{2}(\mathbb{R}^{3})}^{2}+\Vert
\nabla \times \Delta U\Vert _{L^{2}(\mathbb{R}^{3})}^{2}\right) \\
& +\frac{\gamma }{4}\left( \Vert \nabla B\Vert _{L^{2}(\mathbb{R}%
^{3})}^{2}+\Vert \Delta B\Vert _{L^{2}(\mathbb{R}^{3})}^{2}+\Vert \nabla
\times \Delta B\Vert _{L^{2}(\mathbb{R}^{3})}^{2}+\Vert div\,\Delta B\Vert
_{L^{2}(\mathbb{R}^{3})}^{2}\right) \\
& \leq C\left( \Vert U\Vert _{H^{2}}^{2}+\Vert B\Vert _{H^{2}}^{2}\right)
^{3}+C\left( \Vert U\Vert _{H^{2}}^{2}+\Vert B\Vert _{H^{2}}^{2}\right) L(t).
\end{align*}
\end{rem}
\end{lee}

\subsection{\protect\bigskip Proof of Lemma \protect\ref{10}}

\label{proof-lem10}

\hspace{0.6cm}We multiply the first and second equations in (\ref{3}) by $u$
and $b,$ respectively, and afterwards we integrate to obtain
\begin{align}
\frac{1}{2}\frac{d}{dt}\left( \Vert u\Vert _{L^{2}}^{2}+\Vert b\Vert
_{L^{2}}^{2}\right) & +\mu \Vert \nabla u\Vert _{L^{2}}^{2}+\gamma \Vert
\nabla b\Vert _{L^{2}}^{2}=-\overbrace{(P[[u.\nabla ]u],u)_{L^{2}}}^{I_{1}}
\notag \\
& +\overbrace{(P[(\nabla \times b)\times b],u)_{L^{2}}}^{I_{2}}+\overbrace{%
(\nabla \times (u\times b),b)_{L^{2}}}^{I_{3}}  \notag \\
& -\overbrace{(\nabla \times ((\nabla \times b)\times b),b)_{L^{2}}}^{I_{4}}.
\label{aux-eq-1}
\end{align}%
Performing an integration by parts and using the incompressible condition to
$u$, we get
\begin{align}
I_{1}& =([u.\nabla ]u,u)_{L^{2}}=-([u.\nabla ]u,u)_{L^{2}}=0,  \notag \\
&  \notag \\
I_{2}& =((\nabla \times b)\times b,u)_{L^{2}}=([b.\nabla ]b-\frac{1}{2}%
\nabla |b|^{2},u)_{L^{2}}=([b.\nabla ]b,u)_{L^{2}},  \notag \\
&  \notag \\
I_{3}& =(\nabla \times (u\times b),b)_{L^{2}}=(u(div\,b)+[b.\nabla
]u-[u.\nabla ]b,b)_{L^{2}}=(u(div\,b)+[b.\nabla ]u,b)_{L^{2}}  \notag \\
& =-([b.\nabla ]b,u)_{L^{2}}=-I_{2},  \notag \\
&  \notag \\
I_{4}& =((\nabla \times b)\times b,\nabla \times b)_{L^{2}}=0.  \notag
\end{align}%
Inserting the above equalities in (\ref{aux-eq-1}), we obtain (\ref{in1}).

Now, we multiply the first equation of (\ref{3}) by $-\Delta u$ and the
second by $-\Delta b$ in order to obtain
\begin{align}
\frac{1}{2}\frac{d}{dt}\left( \Vert \nabla u\Vert _{L^{2}}^{2}+\Vert \nabla
b\Vert _{L^{2}}^{2}\right) & +\mu \Vert \Delta u\Vert _{L^{2}}^{2}+\gamma
\Vert \Delta b\Vert _{L^{2}}^{2}=\overbrace{(P[[u.\nabla ]u],\Delta
u)_{L^{2}}}^{J_{1}}  \notag \\
& -\overbrace{(P[(\nabla \times b)\times b],\Delta u)_{L^{2}}}^{J_{2}}-%
\overbrace{(\nabla \times (u\times b),\Delta b)_{L^{2}}}^{J_{3}}  \notag \\
& +\overbrace{(\nabla \times ((\nabla \times b)\times b),\Delta b)_{L^{2}}}%
^{J_{4}}.  \label{aux-eq-2}
\end{align}%
The terms $J_{i}$ on the right side of (\ref{aux-eq-2}) can be estimated by
using Gagliardo-Nirenberg inequalities and Young inequality with $\epsilon $%
. Hereafter $C$ will be positive constants that can change in each line and
they may depend of $\mu $, $\gamma $ and $\epsilon $. Precisely, we have%
\begin{align}
|J_{1}|& =|([u.\nabla ]u,\Delta u)_{L^{2}}|\leq C\Vert \Delta u\Vert
_{L^{2}}\Vert \nabla u\Vert _{L^{3}}\Vert u\Vert _{L^{6}}  \notag \\
& \leq C\Vert \Delta u\Vert _{L^{2}}\left( \Vert \nabla u\Vert _{L^{2}}^{%
\frac{1}{2}}\Vert \Delta u\Vert _{L^{2}}^{\frac{1}{2}}\right) \Vert \nabla
u\Vert _{L^{2}}\leq \epsilon \Vert \Delta u\Vert _{L^{2}}^{2}+C\Vert \nabla
u\Vert _{L^{2}}^{6},  \notag \\
&  \notag \\
|J_{2}|& =|((\nabla \times b)\times b,\Delta u)_{L^{2}}|\leq C\Vert \Delta
u\Vert _{L^{2}}\Vert \nabla b\Vert _{L^{3}}\Vert b\Vert _{L^{6}}  \notag \\
& \leq C\Vert \Delta u\Vert _{L^{2}}\left( \Vert \nabla b\Vert _{L^{2}}^{%
\frac{1}{2}}\Vert \Delta b\Vert _{L^{2}}^{\frac{1}{2}}\right) \Vert \nabla
b\Vert _{L^{2}}\leq \epsilon \Vert \Delta u\Vert _{L^{2}}^{2}+\epsilon \Vert
\Delta b\Vert _{L^{2}}^{2}+C\Vert \nabla b\Vert _{L^{2}}^{6},  \notag \\
&  \notag \\
|J_{3}|& =|(\nabla \times (u\times b),\Delta b)_{L^{2}}|\leq C\Vert \Delta
b\Vert _{L^{2}}\left( \Vert \nabla b\Vert _{L^{3}}\Vert u\Vert
_{L^{6}}+\Vert \nabla u\Vert _{L^{2}}\Vert b\Vert _{L^{\infty }}\right)
\notag \\
& \leq C\Vert \Delta b\Vert _{L^{2}}\left( \Vert \nabla b\Vert _{L^{2}}^{%
\frac{1}{2}}\Vert \Delta b\Vert _{L^{2}}^{\frac{1}{2}}\right) \Vert \nabla
u\Vert _{L^{2}}\leq \epsilon \Vert \Delta u\Vert _{L^{2}}^{2}+C\Vert \nabla
u\Vert _{L^{2}}^{4}\Vert \nabla b\Vert _{L^{2}}^{2},  \notag \\
&  \notag \\
|J_{4}|& =|(\nabla \times ((\nabla \times b)\times b),\Delta b)_{L^{2}}|\leq
C\Vert \Delta b\Vert _{L^{2}}\left( \Vert \nabla b\Vert _{L^{3}}\Vert \nabla
b\Vert _{L^{6}}+\Vert b\Vert _{L^{\infty }}\Vert \Delta b\Vert
_{L^{2}}\right)  \notag \\
& \leq C\Vert \Delta b\Vert _{L^{2}}\left( \Vert \nabla b\Vert _{L^{2}}^{%
\frac{1}{2}}\Vert \Delta b\Vert _{L^{2}}^{\frac{1}{2}}\right) \Vert \Delta
b\Vert _{L^{2}}\leq \epsilon \Vert \Delta b\Vert _{L^{2}}^{2}+C\Vert \Delta
b\Vert _{L^{2}}^{4}\Vert \nabla b\Vert _{L^{2}}^{2}.  \notag
\end{align}%
The inequality (\ref{in2}) follows by inserting the above estimates in (\ref%
{aux-eq-2}) and choosing $\epsilon >0$ small enough.

Finally, we deal with (\ref{in4}). After applying $\nabla \times $ in (\ref%
{3}), we multiply the first equation by $-\nabla \times \Delta u$ and the
second by $-\nabla \times \Delta b$. Also we apply $div$ in the second
equation of (\ref{3}) and then multiply it by $-div\,\Delta b.$ With these
manipulations, we get the following equality
\begin{align}
\frac{1}{2}\frac{d}{dt}(\Vert \Delta u\Vert _{L^{2}}^{2}+\Vert \Delta b\Vert
_{L^{2}}^{2})& +\mu \Vert \nabla \times \Delta u\Vert _{L^{2}}^{2}+\gamma
\left( \Vert \nabla \times \Delta b\Vert _{L^{2}}^{2}+\Vert div\,\Delta
b\Vert _{L^{2}}^{2}\right)  \notag \\
& =\overbrace{(\nabla \times ([u.\nabla ]u),\nabla \times \Delta u)_{L^{2}}}%
^{K_{1}}-\overbrace{(\nabla \times ((\nabla \times b)\times b),\nabla \times
\Delta u)_{L^{2}}}^{K_{2}}  \notag \\
& -\overbrace{(\nabla \times \nabla \times (u\times b),\nabla \times \Delta
b)_{L^{2}}}^{K_{3}}+\overbrace{(\nabla \times \nabla \times ((\nabla \times
b)\times b),\nabla \times \Delta b)_{L^{2}}}^{K_{4}}.  \label{aux-eq-3}
\end{align}
For the first three terms on the right side of (\ref{aux-eq-3}), we have
\begin{align}
|K_{1}|& \leq C\Vert \nabla \times \Delta u\Vert _{L^{2}}\left( \Vert u\Vert
_{L^{6}}\Vert \Delta u\Vert _{L^{3}}+\Vert \nabla u\Vert _{L^{\infty }}\Vert
\nabla u\Vert _{L^{2}}\right)  \notag \\
& \leq C\Vert \nabla \times \Delta u\Vert _{L^{2}}\left( \Vert \nabla \times
\Delta u\Vert _{L^{2}}^{\frac{1}{2}}\Vert \Delta u\Vert _{L^{2}}^{\frac{1}{2}%
}\right) \Vert \nabla u\Vert _{L^{2}}  \notag \\
& \leq \epsilon \Vert \nabla \times \Delta u\Vert _{L^{2}}^{2}+C\Vert \Delta
u\Vert _{L^{2}}^{2}\Vert \nabla u\Vert _{L^{2}}^{4},  \notag
\end{align}%
\begin{align}
|K_{2}|& \leq C\Vert \nabla \times \Delta u\Vert _{L^{2}}\left( \Vert b\Vert
_{L^{\infty }}\Vert \Delta b\Vert _{L^{2}}+\Vert \nabla b\Vert _{L^{3}}\Vert
\nabla b\Vert _{L^{6}}\right)  \notag \\
& \leq C\Vert \nabla \times \Delta u\Vert _{L^{2}}\left( \Vert \nabla b\Vert
_{L^{2}}^{\frac{1}{2}}\Vert \Delta b\Vert _{L^{2}}^{\frac{1}{2}}\right)
\Vert \Delta b\Vert _{L^{2}}  \notag \\
& \leq \epsilon \Vert \nabla \times \Delta u\Vert _{L^{2}}^{2}+\epsilon
\Vert \Delta b\Vert _{L^{2}}^{2}+C\Vert \Delta b\Vert _{L^{2}}^{4}\Vert
\nabla b\Vert _{L^{2}}^{2},  \notag
\end{align}%
and
\begin{align}
|K_{3}|& \leq C\Vert \nabla \times \Delta b\Vert _{L^{2}}\left( \Vert u\Vert
_{L^{6}}\Vert \Delta b\Vert _{L^{3}}+\Vert \nabla u\Vert _{L^{2}}\Vert
\nabla b\Vert _{L^{\infty }}\right)  \notag \\
& +C\Vert \nabla \times \Delta b\Vert _{L^{2}}\left( \Vert b\Vert
_{L^{6}}\Vert \Delta u\Vert _{L^{3}}\right)  \notag \\
& \leq C\Vert \nabla \times \Delta b\Vert _{L^{2}}\left( \Vert \Delta b\Vert
_{L^{2}}^{\frac{1}{2}}\left( \Vert \nabla \times \Delta b\Vert _{L^{2}}^{%
\frac{1}{2}}+\Vert div\,\Delta b\Vert _{L^{2}}^{\frac{1}{2}}\right) \right)
\Vert \nabla u\Vert _{L^{2}}  \notag \\
& +C\Vert \nabla \times \Delta b\Vert _{L^{2}}\left( \Vert \nabla \times
\Delta u\Vert _{L^{2}}^{\frac{1}{2}}\Vert \Delta u\Vert _{L^{2}}^{\frac{1}{2}%
}\right) \Vert \nabla b\Vert _{L^{2}}  \notag \\
& \leq \epsilon \Vert \nabla \times \Delta u\Vert _{L^{2}}^{2}+\epsilon
\Vert \nabla \times \Delta b\Vert _{L^{2}}^{2}+\epsilon \Vert div\,\Delta
b\Vert _{L^{2}}^{2}+C\Vert \Delta b\Vert _{L^{2}}^{2}\Vert \nabla u\Vert
_{L^{2}}^{4}+C\Vert \Delta u\Vert _{L^{2}}^{2}\Vert \nabla b\Vert
_{L^{2}}^{4}.  \notag
\end{align}%
For the Hall-term, using the vector identity $(A\times B).A=0$, we get
\begin{align}
K_{4}& =(\nabla \times \nabla \times ((\nabla \times b)\times b),-\nabla
\times \nabla \times \nabla \times b)_{L^{2}}  \notag \\
& =\overbrace{(\nabla \times \nabla \times ((\nabla \times b)\times
b)-\nabla \times ((\nabla \times \nabla \times b)\times b),-\nabla \times
\nabla \times \nabla \times b)_{L^{2}}}^{K_{5}}  \notag \\
& +\overbrace{(\nabla \times ((\nabla \times \nabla \times b)\times
b)-(\nabla \times \nabla \times \nabla \times b)\times b,-\nabla \times
\nabla \times \nabla \times b)_{L^{2}}}^{K_{6}}.  \label{aux-Hall-1}
\end{align}%
Now, by identities (\ref{345-1})-(\ref{345-2}), we obtain
\begin{align}
|K_{5}|& =|(\nabla \times \{\nabla \times ((\nabla \times b)\times
b)-(\nabla \times \nabla \times b)\times b\},\nabla \times \Delta b)_{L^{2}}|
\notag \\
& =|(\nabla \times \{(\nabla \times b)(div\,b)-2[(\nabla \times b).\nabla
]b\},\nabla \times \Delta b)_{L^{2}}|  \notag \\
& \leq C\Vert \nabla \times \Delta b\Vert _{L^{2}}\Vert \nabla b\Vert
_{L^{6}}\Vert \Delta b\Vert _{L^{3}}  \notag \\
& \leq C\Vert \nabla \times \Delta b\Vert _{L^{2}}\Vert \Delta b\Vert
_{L^{2}}\left( \Vert \Delta b\Vert _{L^{2}}^{\frac{1}{2}}\left( \Vert \nabla
\times \Delta b\Vert _{L^{2}}^{\frac{1}{2}}+\Vert div\,\Delta b\Vert
_{L^{2}}^{\frac{1}{2}}\right) \right)  \notag \\
& \leq \epsilon \Vert \nabla \times \Delta b\Vert _{L^{2}}^{2}+\epsilon
\Vert div\,\Delta b\Vert _{L^{2}}^{2}+C\Vert \Delta b\Vert _{L^{2}}^{6}
\label{aux-Hall-2}
\end{align}%
and%
\begin{align}
|K_{6}|& =|(\nabla \times \nabla \times b)(div\,b)-2[(\nabla \times \nabla
\times b).\nabla ]b-(\nabla \times \nabla \times b)\times (\nabla \times
b),\nabla \times \Delta b)_{L^{2}}|  \notag \\
& \leq C\Vert \nabla \times \Delta b\Vert _{L^{2}}\Vert \Delta b\Vert
_{L^{3}}\Vert \nabla b\Vert _{L^{6}}  \notag \\
& \leq C\Vert \nabla \times \Delta b\Vert _{L^{2}}\left( \Vert \Delta b\Vert
_{L^{2}}^{\frac{1}{2}}\left( \Vert \nabla \times \Delta b\Vert _{L^{2}}^{%
\frac{1}{2}}+\Vert div\,\Delta b\Vert _{L^{2}}^{\frac{1}{2}}\right) \right)
\Vert \Delta b\Vert _{L^{2}}  \notag \\
& \leq \epsilon \Vert \nabla \times \Delta b\Vert _{L^{2}}^{2}+\epsilon
\Vert div\,\Delta b\Vert _{L^{2}}^{2}+C\Vert \Delta b\Vert _{L^{2}}^{6}.
\label{aux-Hall-3}
\end{align}%
We obtain (\ref{in4}) from (\ref{aux-eq-3}), the above estimates for $%
\left\vert K_{i}\right\vert ,$ and choosing a suitable $\epsilon >0$ small
enough.

\begin{flushright}
\rule{2mm}{2mm}
\end{flushright}

%%%%%%%%%%%%%%%%%%%%%%%%%%%%%%%%%%%%%%%%%%%%%%%%%%%%%%%%%%%%%%%%%%%%%%%%%%%%%%%%%%%%%%%%%%%%%%%%%%%%%%%%%

\subsection{\protect\bigskip Proof of Lemma \protect\ref{11}}

\label{proof-lem11}

\hspace{0.6cm}The proof of (\ref{4}) is a straightforward calculation. Now
we multiply the first equation of (\ref{4}) by $U$ and the second equation
by $B$ to obtain
\begin{align}
\frac{1}{2}\frac{d}{dt}\left( \Vert U\Vert _{L^{2}}^{2}\right. & +\left.
\Vert B\Vert _{L^{2}}^{2}\right) +\mu \Vert \nabla U\Vert
_{L^{2}}^{2}+\gamma \Vert \nabla B\Vert _{L^{2}}^{2}=-\overbrace{%
(P[[U.\nabla ]U],U)_{L^{2}}}^{L_{1}}-\overbrace{(P[[U.\nabla ]u],U)_{L^{2}}}%
^{L_{2}}  \notag \\
& -\overbrace{(P[[u.\nabla ]U],U)_{L^{2}}}^{L_{3}}+\overbrace{(P[(\nabla
\times B)\times B],U)_{L^{2}}}^{L_{4}}+\overbrace{(P[(\nabla \times B)\times
b],U)_{L^{2}}}^{L_{5}}  \notag \\
& +\overbrace{(P[(\nabla \times b)\times B],U)_{L^{2}}}^{L_{6}}+\overbrace{%
(\nabla \times (U\times B),B)_{L^{2}}}^{L_{7}}+\overbrace{(\nabla \times
(U\times b),B)_{L^{2}}}^{L_{8}}  \notag \\
& +\overbrace{(\nabla \times (u\times B),B)_{L^{2}}}^{L_{9}}-\overbrace{%
(\nabla \times ((\nabla \times B)\times B),B)_{L^{2}}}^{L_{10}}-\overbrace{%
(\nabla \times ((\nabla \times B)\times b),B)_{L^{2}}}^{L_{11}}  \notag \\
& -\overbrace{(\nabla \times ((\nabla \times b)\times B),B)_{L^{2}}}%
^{L_{12}}.  \label{aux-eq-5}
\end{align}%
We have that $L_{1}=L_{3}=L_{10}=L_{11}=0$ and $L_{4}=-L_{7}.$ Then, we need
to estimate the remainder terms. Proceeding as in the proof of Lemma \ref{10}%
, we obtain
\begin{align}
|L_{2}|& \leq \epsilon \Vert \nabla U\Vert _{L^{2}}^{2}+C\Vert \nabla u\Vert
_{L^{2}}^{4}\Vert U\Vert _{L^{2}}^{2},  \notag \\
|L_{5}|& \leq \epsilon \Vert \nabla B\Vert _{L^{2}}^{2}+C\Vert \nabla b\Vert
_{L^{2}}^{4}\Vert B\Vert _{L^{2}}^{2},  \notag \\
|L_{6}|& \leq \epsilon \Vert \nabla B\Vert _{L^{2}}^{2}+\epsilon \Vert
\nabla U\Vert _{L^{2}}^{2}+C\Vert \nabla b\Vert _{L^{2}}^{4}\Vert B\Vert
_{L^{2}}^{2},  \notag \\
|L_{8}|& \leq \epsilon \Vert \nabla B\Vert _{L^{2}}^{2}+\epsilon \Vert
\nabla U\Vert _{L^{2}}^{2}+C\Vert \nabla b\Vert _{L^{2}}^{4}\Vert B\Vert
_{L^{2}}^{2},  \notag \\
|L_{9}|& \leq \epsilon \Vert \nabla B\Vert _{L^{2}}^{2}+C\Vert \nabla u\Vert
_{L^{2}}^{4}\Vert B\Vert _{L^{2}}^{2},  \notag \\
|L_{12}|& \leq \epsilon \Vert \nabla B\Vert _{L^{2}}^{2}+C\Vert \Delta
b\Vert _{L^{2}}^{4}\Vert B\Vert _{L^{2}}^{2}.  \notag
\end{align}%
Now we obtain (\ref{in5}) after inserting the above inequalities in (\ref%
{aux-eq-5}) and taking $\epsilon >0$ small enough.

Next we prove (\ref{in50}). Multiplying the first equation of (\ref{4}) by $%
-\Delta U$ and the second by $-\Delta B,$ we get

\begin{align}
\frac{1}{2}\frac{d}{dt}\left( \Vert \nabla U\Vert _{L^{2}}^{2}\right. &
+\left. \Vert \nabla B\Vert _{L^{2}}^{2}\right) +\mu \Vert \Delta U\Vert
_{L^{2}}^{2}+\gamma \Vert \Delta B\Vert _{L^{2}}^{2}=\overbrace{(P[[U.\nabla
]U],\Delta U)_{L^{2}}}^{M_{1}}+\overbrace{(P[[U.\nabla ]u],\Delta U)_{L^{2}}}%
^{M_{2}}  \notag \\
& +\overbrace{(P[[u.\nabla ]U],\Delta U)_{L^{2}}}^{M_{3}}-\overbrace{%
(P[(\nabla \times B)\times B],\Delta U)_{L^{2}}}^{M_{4}}-\overbrace{%
(P[(\nabla \times B)\times b],\Delta U)_{L^{2}}}^{M_{5}}  \notag \\
& -\overbrace{(P[(\nabla \times b)\times B],\Delta U)_{L^{2}}}^{M_{6}}-%
\overbrace{(\nabla \times (U\times B),\Delta B)_{L^{2}}}^{M_{7}}-\overbrace{%
(\nabla \times (U\times b),\Delta B)_{L^{2}}}^{M_{8}}  \notag \\
& -\overbrace{(\nabla \times (u\times B),\Delta B)_{L^{2}}}^{M_{9}}+%
\overbrace{(\nabla \times ((\nabla \times B)\times B),\Delta B)_{L^{2}}}%
^{M_{10}}  \notag \\
& +\overbrace{(\nabla \times ((\nabla \times B)\times b),\Delta B)_{L^{2}}}%
^{M_{11}}+\overbrace{(\nabla \times ((\nabla \times b)\times B),\Delta
B)_{L^{2}}}^{M_{12}}.  \label{aux-eq-6}
\end{align}

We have the following estimates:
\begin{eqnarray*}
|M_{1}| &\leq &\epsilon \Vert \Delta U\Vert _{L^{2}}^{2}+C\Vert \nabla
U\Vert _{L^{2}}^{6}, \\
|M_{2}|&\leq& \epsilon \Vert \Delta U\Vert _{L^{2}}^{2}+C\Vert \nabla u\Vert
_{L^{2}}^{4}\Vert \nabla U\Vert _{L^{2}}^{2}, \\
|M_{3}| &\leq &\epsilon \Vert \Delta U\Vert _{L^{2}}^{2}+C\Vert \nabla
u\Vert _{L^{2}}^{4}\Vert \nabla U\Vert _{L^{2}}^{2}, \\
|M_{4}|&\leq& \epsilon \Vert \Delta U\Vert _{L^{2}}^{2}+\epsilon \Vert
\Delta B\Vert _{L^{2}}^{2}+C\Vert \nabla B\Vert _{L^{2}}^{6}, \\
|M_{5}| &\leq &\epsilon \Vert \Delta U\Vert _{L^{2}}^{2}+\epsilon \Vert
\Delta B\Vert _{L^{2}}^{2}+C\Vert \nabla b\Vert _{L^{2}}^{4}\Vert \nabla
B\Vert _{L^{2}}^{2}, \\
|M_{6}|&\leq& \epsilon \Vert \Delta U\Vert _{L^{2}}^{2}+\epsilon \Vert
\Delta B\Vert _{L^{2}}^{2}+C\Vert \nabla b\Vert _{L^{2}}^{4}\Vert \nabla
B\Vert _{L^{2}}^{2}, \\
|M_{7}| &\leq &\epsilon \Vert \Delta U\Vert _{L^{2}}^{2}+\epsilon \Vert
\Delta B\Vert _{L^{2}}^{2}+C\Vert \nabla B\Vert _{L^{2}}^{4}\Vert \nabla
U\Vert _{L^{2}}^{2} \\
|M_{8}|&\leq& \epsilon \Vert \Delta U\Vert _{L^{2}}^{2}+\epsilon \Vert
\Delta B\Vert _{L^{2}}^{2}+C\Vert \nabla b\Vert _{L^{2}}^{4}\Vert \nabla
U\Vert _{L^{2}}^{2}, \\
|M_{9}| &\leq &\epsilon \Vert \Delta B\Vert _{L^{2}}^{2}+C\Vert \nabla
u\Vert _{L^{2}}^{4}\Vert \nabla B\Vert _{L^{2}}^{2}, \\
|M_{10}|&\leq& \epsilon \Vert \Delta B\Vert _{L^{2}}^{2}+C\Vert \Delta
B\Vert _{L^{2}}^{4}\Vert \nabla B\Vert _{L^{2}}^{2},
\end{eqnarray*}%
and
\begin{eqnarray}
|M_{11}| &\leq &\epsilon \Vert \Delta B\Vert _{L^{2}}^{2}+C\Vert \Delta
b\Vert _{L^{2}}^{4}\Vert \nabla B\Vert _{L^{2}}^{2}+C\Vert \Delta B\Vert
_{L^{2}}^{2}\Vert \nabla b\Vert _{L^{2}}^{4}+C\Vert \Delta B\Vert
_{L^{2}}^{2}\Vert \Delta b\Vert _{L^{2}}^{4},  \notag \\
|M_{12}| &\leq &\epsilon \Vert \Delta B\Vert _{L^{2}}^{2}+C\Vert \Delta
b\Vert _{L^{2}}^{4}\Vert \nabla B\Vert _{L^{2}}^{2}.  \notag
\end{eqnarray}%
Again, taking $\epsilon >0$ small enough, the above estimates together with (%
\ref{aux-eq-6}) give (\ref{in50}). \bigskip

\begin{flushright}
  \rule{2mm}{2mm}
 \end{flushright}

\subsection{\protect\bigskip Proof of Lemma \protect\ref{11-b}}

\label{proof-lem11-b}

\hspace{0.6cm}First we apply $\nabla \times $ in (\ref{4}) and then multiply
the first equation by $-\nabla \times \Delta U$ and the second equation by $%
-\nabla \times \Delta B.$ Also, we apply $div$ in the second equation of (%
\ref{4}) and we multiply it by $-div\,\Delta B.$ After this manipulations,
we obtain
\begin{eqnarray}
\frac{1}{2}\frac{d}{dt}\left( \Vert \Delta U\Vert _{L^{2}}^{2}\right.
&+&\left. \Vert \Delta B\Vert _{L^{2}}^{2}\right) +\mu \Vert \nabla \times
\Delta U\Vert _{L^{2}}^{2}+\gamma \left( \Vert \nabla \times \Delta B\Vert
_{L^{2}}^{2}+\Vert div\,\Delta B\Vert _{L^{2}}^{2}\right)  \notag \\
&=&\overbrace{(\nabla \times ([U.\nabla ]U),\nabla \times \Delta U)_{L^{2}}}%
^{N_{1}}+\overbrace{(\nabla \times ([U.\nabla ]u),\nabla \times \Delta
U)_{L^{2}}}^{N_{2}}  \notag \\
&+&\overbrace{(\nabla \times ([u.\nabla ]U),\nabla \times \Delta U)_{L^{2}}}%
^{N_{3}}-\overbrace{(\nabla \times ((\nabla \times B)\times B),\nabla \times
\Delta U)_{L^{2}}}^{N_{4}}  \notag \\
&-&\overbrace{(\nabla \times ((\nabla \times B)\times b),\nabla \times
\Delta U)_{L^{2}}}^{N_{5}}-\overbrace{(\nabla \times ((\nabla \times
b)\times B),\nabla \times \Delta U)_{L^{2}}}^{N_{6}}  \notag \\
&-&\overbrace{(\nabla \times \nabla \times (U\times B),\nabla \times \Delta
B)_{L^{2}}}^{N_{7}}-\overbrace{(\nabla \times \nabla \times (U\times
b),\nabla \times \Delta B)_{L^{2}}}^{N_{8}}  \notag \\
&-&\overbrace{(\nabla \times \nabla \times (u\times B),\nabla \times \Delta
B)_{L^{2}}}^{N_{9}}+\overbrace{(\nabla \times \nabla \times ((\nabla \times
B)\times B),\nabla \times \Delta B)_{L^{2}}}^{N_{10}}  \notag \\
&+&\overbrace{(\nabla \times \nabla \times ((\nabla \times B)\times
b),\nabla \times \Delta B)_{L^{2}}}^{N_{11}}+\overbrace{(\nabla \times
\nabla \times ((\nabla \times b)\times B),\nabla \times \Delta B)_{L^{2}}}%
^{N_{12}}.  \notag \\
&&  \label{aux-Hall-1-1}
\end{eqnarray}%
As before, we can estimate
\begin{eqnarray}
|N_{1}| &\leq &\epsilon \Vert \nabla \times \Delta U\Vert
_{L^{2}}^{2}+\epsilon \Vert \Delta U\Vert _{L^{2}}^{2}+C\Vert \Delta U\Vert
_{L^{2}}^{4}\Vert \nabla U\Vert _{L^{2}}^{2},  \label{aux-Hall-1-2} \\
|N_{2}| &\leq &\epsilon \Vert \nabla \times \Delta U\Vert
_{L^{2}}^{2}+\epsilon \Vert \Delta U\Vert _{L^{2}}^{2}+C\Vert \Delta u\Vert
_{L^{2}}^{4}\Vert \nabla U\Vert _{L^{2}}^{2},  \notag \\
|N_{3}| &\leq &\epsilon \Vert \nabla \times \Delta U\Vert
_{L^{2}}^{2}+C\Vert \nabla u\Vert _{L^{2}}^{4}\Vert \Delta U\Vert
_{L^{2}}^{2},  \notag \\
|N_{4}| &\leq &\epsilon \Vert \nabla \times \Delta U\Vert
_{L^{2}}^{2}+\epsilon \Vert \Delta B\Vert _{L^{2}}^{2}+C\Vert \nabla B\Vert
_{L^{2}}^{2}\Vert \Delta B\Vert _{L^{2}}^{4},  \notag \\
|N_{5}| &\leq &\epsilon \Vert \nabla \times \Delta U\Vert
_{L^{2}}^{2}+\epsilon \Vert \nabla \times \Delta B\Vert
_{L^{2}}^{2}+\epsilon \Vert div\,\Delta B\Vert _{L^{2}}^{2}+C\Vert \nabla
b\Vert _{L^{2}}^{4}\Vert \Delta B\Vert _{L^{2}}^{2},  \notag
\end{eqnarray}%
and%
\begin{eqnarray}
|N_{6}| &\leq &\epsilon \Vert \nabla \times \Delta U\Vert
_{L^{2}}^{2}+\epsilon \Vert \Delta B\Vert _{L^{2}}^{2}+C\Vert \Delta b\Vert
_{L^{2}}^{4}\Vert \nabla B\Vert _{L^{2}}^{2},  \label{aux-Hall-1-3} \\
|N_{7}| &\leq &\epsilon \Vert \Delta U\Vert _{L^{2}}^{2}+\epsilon \Vert
\Delta B\Vert _{L^{2}}^{2}+\epsilon \Vert \nabla \times \Delta B\Vert
_{L^{2}}^{2}+C\Vert \Delta B\Vert _{L^{2}}^{4}\Vert \nabla U\Vert
_{L^{2}}^{2}+C\Vert \Delta U\Vert _{L^{2}}^{4}\Vert \nabla B\Vert
_{L^{2}}^{2},  \notag \\
|N_{8}| &\leq &\epsilon \Vert \nabla \times \Delta B\Vert
_{L^{2}}^{2}+\epsilon \Vert \nabla \times \Delta U\Vert
_{L^{2}}^{2}+\epsilon \Vert \Delta U\Vert _{L^{2}}^{2}+C\Vert \nabla b\Vert
_{L^{2}}^{4}\Vert \Delta U\Vert _{L^{2}}^{2}+C\Vert \nabla U\Vert
_{L^{2}}^{2}\Vert \Delta b\Vert _{L^{2}}^{4},  \notag \\
|N_{9}| &\leq &\epsilon \Vert \nabla \times \Delta B\Vert
_{L^{2}}^{2}+\epsilon \Vert div\,\Delta B\Vert _{L^{2}}^{2}+\epsilon \Vert
\Delta B\Vert _{L^{2}}^{2}+C\Vert \nabla u\Vert _{L^{2}}^{4}\Vert \Delta
B\Vert _{L^{2}}^{2}+C\Vert \nabla B\Vert _{L^{2}}^{2}\Vert \Delta u\Vert
_{L^{2}}^{4}.  \notag
\end{eqnarray}%
For the Hall-term, using the vector identity $(A\times B).A=0$, we get
\begin{eqnarray}
N_{10} &=&(\nabla \times \nabla \times ((\nabla \times B)\times B),-\nabla
\times \nabla \times \nabla \times B)_{L^{2}}  \notag \\
&=&\overbrace{(\nabla \times \nabla \times ((\nabla \times B)\times
B)-\nabla \times ((\nabla \times \nabla \times B)\times B),-\nabla \times
\nabla \times \nabla \times B)_{L^{2}}}^{N_{10,a}}  \notag \\
&+&\overbrace{(\nabla \times ((\nabla \times \nabla \times B)\times
B)-(\nabla \times \nabla \times \nabla \times B)\times B,-\nabla \times
\nabla \times \nabla \times B)_{L^{2}}}^{N_{10,b}}  \label{aux-Hall-4}
\end{eqnarray}%
and%
\begin{eqnarray}
N_{11} &=&(\nabla \times \nabla \times ((\nabla \times B)\times b),-\nabla
\times \nabla \times \nabla \times B)_{L^{2}}  \notag \\
&=&\overbrace{(\nabla \times \nabla \times ((\nabla \times B)\times
b)-\nabla \times ((\nabla \times \nabla \times B)\times b),-\nabla \times
\nabla \times \nabla \times B)_{L^{2}}}^{N_{11,a}}  \notag \\
&+&\overbrace{(\nabla \times ((\nabla \times \nabla \times B)\times
b)-(\nabla \times \nabla \times \nabla \times B)\times b,-\nabla \times
\nabla \times \nabla \times B)_{L^{2}}}^{N_{11,b}}.  \label{aux-Hall-5}
\end{eqnarray}%
Now, by using identities (\ref{345-1})-(\ref{345-2}), we obtain
\begin{eqnarray}
|N_{10,a}| &=&|(\nabla \times \{\nabla \times ((\nabla \times B)\times
B)-(\nabla \times \nabla \times B)\times B\},\nabla \times \Delta B)_{L^{2}}|
\notag \\
&=&|(\nabla \times \{(\nabla \times B)(div\,B)-2[(\nabla \times B).\nabla
]B\},\nabla \times \Delta B)_{L^{2}}|  \notag \\
&\leq &C\Vert \Delta B\Vert _{L^{2}}\Vert \nabla B\Vert _{L^{\infty }}\Vert
\nabla \times \Delta B\Vert _{L^{2}}  \notag \\
&\leq &C\Vert \Delta B\Vert _{L^{2}}^{\frac{3}{2}}\left( \Vert \nabla \times
\Delta B\Vert _{L^{2}}^{\frac{1}{2}}+\Vert div\,\Delta B\Vert _{L^{2}}^{%
\frac{1}{2}}\right) \Vert \nabla \times \Delta B\Vert _{L^{2}}  \notag \\
&\leq &\epsilon \Vert \nabla \times \Delta B\Vert _{L^{2}}^{2}+\epsilon
\Vert div\,\Delta B\Vert _{L^{2}}^{2}+C\Vert \Delta B\Vert _{L^{2}}^{6},
\label{aux-Hall-6}
\end{eqnarray}%
\begin{eqnarray}
|N_{10,b}| &=&|((\nabla \times \nabla \times B)(div\,B)-2[(\nabla \times
\nabla \times B).\nabla ]B-(\nabla \times \nabla \times B)\times (\nabla
\times B),\nabla \times \Delta B)_{L^{2}}|  \notag \\
&\leq &C\Vert \Delta B\Vert _{L^{2}}\Vert \nabla B\Vert _{L^{\infty }}\Vert
\nabla \times \Delta B\Vert _{L^{2}}  \notag \\
&\leq &C\Vert \Delta B\Vert _{L^{2}}^{\frac{3}{2}}\left( \Vert \nabla \times
\Delta B\Vert _{L^{2}}^{\frac{1}{2}}+\Vert div\,\Delta B\Vert _{L^{2}}^{%
\frac{1}{2}}\right) \Vert \nabla \times \Delta B\Vert _{L^{2}}  \notag \\
&\leq &\epsilon \Vert \nabla \times \Delta B\Vert _{L^{2}}^{2}+\epsilon
\Vert div\,\Delta B\Vert _{L^{2}}^{2}+C\Vert \Delta B\Vert _{L^{2}}^{6},
\label{aux-Hall-7}
\end{eqnarray}%
\begin{eqnarray}
|N_{11,a}| &=&|(\nabla \times \{(\nabla \times B)(div\,b)-2[(\nabla \times
B).\nabla ]b-(\nabla \times B)\times (\nabla \times b)\},\nabla \times
\Delta B)_{L^{2}}|  \notag \\
&\leq &C\left( \Vert \Delta B\Vert _{L^{3}}\Vert \nabla b\Vert
_{L^{6}}+\Vert \Delta b\Vert _{L^{2}}\Vert \nabla B\Vert _{L^{\infty
}}\right) \Vert \nabla \times \Delta B\Vert _{L^{2}}  \notag \\
&\leq &C\Vert \Delta b\Vert _{L^{2}}\Vert \Delta B\Vert _{L^{2}}^{\frac{1}{2}%
}\left( \Vert \nabla \times \Delta B\Vert _{L^{2}}^{\frac{1}{2}}+\Vert
div\,\Delta B\Vert _{L^{2}}^{\frac{1}{2}}\right) \Vert \nabla \times \Delta
B\Vert _{L^{2}}  \notag \\
&\leq &\epsilon \Vert \nabla \times \Delta B\Vert _{L^{2}}^{2}+\epsilon
\Vert div\,\Delta B\Vert _{L^{2}}^{2}+C\Vert \Delta B\Vert _{L^{2}}^{2}\Vert
\Delta b\Vert _{L^{2}}^{4},  \label{aux-Hall-8}
\end{eqnarray}%
\begin{eqnarray}
|N_{11,b}| &=&|((\nabla \times \nabla \times B)(div\,b)-2[(\nabla \times
\nabla \times B).\nabla ]b-(\nabla \times \nabla \times B)\times (\nabla
\times b),\nabla \times \Delta B)_{L^{2}}|  \notag \\
&\leq &C\Vert \Delta B\Vert _{L^{3}}\Vert \nabla b\Vert _{L^{6}}\Vert \nabla
\times \Delta B\Vert _{L^{2}}  \notag \\
&\leq &C\Vert \Delta B\Vert _{L^{2}}^{\frac{1}{2}}\left( \Vert \nabla \times
\Delta B\Vert _{L^{2}}^{\frac{1}{2}}+\Vert div\,\Delta B\Vert _{L^{2}}^{%
\frac{1}{2}}\right) \Vert \Delta b\Vert _{L^{2}}\Vert \nabla \times \Delta
B\Vert _{L^{2}}  \notag \\
&\leq &\epsilon \Vert \nabla \times \Delta B\Vert _{L^{2}}^{2}+\epsilon
\Vert div\,\Delta B\Vert _{L^{2}}^{2}+C\Vert \Delta B\Vert _{L^{2}}^{2}\Vert
\Delta b\Vert _{L^{2}}^{4},  \label{aux-Hall-9}
\end{eqnarray}%
and%
\begin{eqnarray}
|N_{12}| &=&|(\nabla \times \nabla \times ((\nabla \times b)\times B),\nabla
\times \Delta B)_{L^{2}}|  \notag \\
&\leq &C\left( \left( \Vert \nabla \times \Delta b\Vert _{L^{2}}+\Vert
div\,\Delta b\Vert _{L^{2}}\right) \Vert B\Vert _{L^{\infty }}+\Vert \Delta
b\Vert _{L^{6}}\Vert \nabla B\Vert _{L^{3}}\right) \Vert \nabla \times
\Delta B\Vert _{L^{2}}  \notag \\
&+&\Vert \Delta B\Vert _{L^{3}}\Vert \nabla b\Vert _{L^{6}}\Vert \nabla
\times \Delta B\Vert _{L^{2}}  \notag \\
&\leq &C\left( \left( \Vert \nabla \times \Delta b\Vert _{L^{2}}+\Vert
div\,\Delta b\Vert _{L^{2}}\right) \Vert \nabla B\Vert _{L^{2}}^{\frac{1}{2}%
}\Vert \Delta B\Vert _{L^{2}}^{\frac{1}{2}}\right) \Vert \nabla \times
\Delta B\Vert _{L^{2}}  \notag \\
&+&\Vert \Delta B\Vert _{L^{2}}^{\frac{1}{2}}\left( \Vert \nabla \times
\Delta B\Vert _{L^{2}}^{\frac{1}{2}}+\Vert div\,\Delta B\Vert _{L^{2}}^{%
\frac{1}{2}}\right) \Vert \Delta b\Vert _{L^{2}}\Vert \nabla \times \Delta
B\Vert _{L^{2}}  \notag \\
&\leq &\epsilon \Vert \nabla \times \Delta B\Vert _{L^{2}}^{2}+\epsilon
\Vert div\,\Delta B\Vert _{L^{2}}^{2}+C\Vert \Delta B\Vert _{L^{2}}^{2}\Vert
\Delta b\Vert _{L^{2}}^{4}  \notag \\
&+&C\left( \Vert \nabla B\Vert _{L^{2}}^{2}+\Vert \Delta B\Vert
_{L^{2}}^{2}\right) \left( \Vert \nabla \times \Delta b\Vert
_{L^{2}}^{2}+\Vert div\,\Delta b\Vert _{L^{2}}^{2}\right) .
\label{aux-Hall-10}
\end{eqnarray}

We conclude the proof of (\ref{in300}) by considering the estimates (\ref%
{aux-Hall-1-2})-(\ref{aux-Hall-10}) in (\ref{aux-Hall-1-1}) and taking a
suitable $\epsilon >0$.
\begin{flushright}
  \rule{2mm}{2mm}
 \end{flushright}

%%%%%%%%%%%%%%%%%%%%%%%%%%%%%%%%%%%%%%%%%%%%%%%%%%%%%%%%%%%%%%%%%%%%%%%%%%%%%%%%%%%%%%%%%%%%%%%%%%%%%%%%%
%%%%%%%%%%%%%%%%%%%%%%%%%%%%%%%%%%%%%%%%%%%%%%%%%%%%%%%%%%%%%%%%%%%%%%%%%%%%%%%%%%%%%%%%%%%%%%%%%%%%%%%%%

\section{\protect\bigskip Proof of Results}

\label{AAA3}

%%%%%%%%%%%%%%%%%%%%%%%%%%%%%%%%%%%%%%%%%%%%%%%%%%%%%%%%%%%%%%%%%%%%%%%%%%%%%%%%%%%%%%%%%%%%%%%%%%%%%%%%%

\subsection{\protect\bigskip Proof of Theorem \protect\ref{a1}}

\hspace{0.6cm}Let $\rho \in C_{0}^{\infty }(\mathbb{R}^{3})$ be a
non-negative radial scalar function with unit integral and $0\leq \rho \leq
1 $ and let $\mathcal{J}_{n}\equiv \rho _{n}\ast $ be the standard
mollifier, where $\rho _{n}(x)=n^{3}\rho (nx)$ and $n\in \mathbb{N}$ (see
\cite{Adams}). We consider the regularized system
\begin{equation}
\left\{
\begin{array}{rcl}
\displaystyle\partial _{t}u+\mathcal{J}_{n}\mathbb{P}[[\mathcal{J}%
_{n}u.\nabla ]\mathcal{J}_{n}u-(\nabla \times \mathcal{J}_{n}b)\times
\mathcal{J}_{n}b] & = & \mu \Delta \mathcal{J}_{n}^{2}u; \\
\displaystyle\partial _{t}b-\nabla \times \mathcal{J}_{n}(\mathcal{J}%
_{n}u\times \mathcal{J}_{n}b)+\nabla \times \mathcal{J}_{n}((\nabla \times
\mathcal{J}_{n}b)\times \mathcal{J}_{n}b) & = & \gamma \Delta \mathcal{J}%
_{n}^{2}b; \\
\mathbb{P}[u] & = & u; \\
(u_{n,\,0},b_{n,\,0}) & = & (\mathcal{J}_{n}u_{0},\mathcal{J}_{n}b_{0}). \\
&  &
\end{array}%
\right.  \label{10z}
\end{equation}
In \cite[Proposition 3.1 ]{DC1} (see also \cite{Bertozzi}), global-in-time
existence of smooth solutions for (\ref{10z}) is obtained for each fixed $%
n\in \mathbb{N}$. More precisely, the above system is considered as an
autonomous infinite-dimensional ODE system in $H_{\sigma }^{m}\times H^{m}$
and it is used the generalized Picard theorem in Banach spaces. The key
point is the suitable way that the original system is mollified by $\mathcal{%
J}_{n}$ which implies simpler estimates for (\ref{10z}). For further details
see \cite{DC1} and \cite{Bertozzi}.

Let $(u_{n},b_{n})$ be a global-in-time smooth solution of (\ref{10z}) with
initial data $(u_{n,\,0},b_{n,\,0})$. Adapting for (\ref{10z}) the \textit{a
priori} estimates contained in Lemma \ref{10} and Remark \ref{8}, we obtain
a constant $C>0$ (independent of $n$) such that

\begin{align}
& \frac{1}{2}\frac{d}{dt}\left( \Vert u_{n}\Vert _{H^{2}(\mathbb{R}%
^{3})}^{2}+\Vert b_{n}\Vert _{H^{2}(\mathbb{R}^{3})}^{2}\right) +\frac{\mu }{%
4}\left( \Vert \nabla \mathcal{J}_{n}u_{n}\Vert _{L^{2}(\mathbb{R}%
^{3})}^{2}+\Vert \Delta \mathcal{J}_{n}u_{n}\Vert _{L^{2}(\mathbb{R}%
^{3})}^{2}+\Vert \nabla \times \Delta \mathcal{J}_{n}u_{n}\Vert _{L^{2}(%
\mathbb{R}^{3})}^{2}\right)  \notag \\
& +\frac{\gamma }{4}\left( \Vert \nabla \mathcal{J}_{n}b_{n}\Vert _{L^{2}(%
\mathbb{R}^{3})}^{2}+\Vert \Delta \mathcal{J}_{n}b_{n}\Vert _{L^{2}(\mathbb{R%
}^{3})}^{2}+\Vert \nabla \times \Delta \mathcal{J}_{n}b_{n}\Vert _{L^{2}(%
\mathbb{R}^{3})}^{2}+\Vert div\,\Delta \mathcal{J}_{n}b_{n}\Vert _{L^{2}(%
\mathbb{R}^{3})}^{2}\right)  \notag \\
& \leq C\left( \Vert u_{n}\Vert _{H^{2}(\mathbb{R}^{3})}^{2}+\Vert
b_{n}\Vert _{H^{2}(\mathbb{R}^{3})}^{2}\right) \left( \Vert \nabla
u_{n}\Vert _{L^{2}(\mathbb{R}^{3})}^{4}+\Vert \nabla b_{n}\Vert _{L^{2}(%
\mathbb{R}^{3})}^{4}+\Vert \Delta b_{n}\Vert _{L^{2}(\mathbb{R}%
^{3})}^{4}\right) .  \label{cach2n}
\end{align}%
It follows that
\begin{equation*}
\frac{1}{2}\frac{d}{dt}\left( \Vert u_{n}(t)\Vert _{H^{2}(\mathbb{R}%
^{3})}^{2}+\Vert b_{n}(t)\Vert _{H^{2}(\mathbb{R}^{3})}^{2}\right) \leq
C\left( \Vert u_{n}(t)\Vert _{H^{2}(\mathbb{R}^{3})}^{2}+\Vert b_{n}(t)\Vert
_{H^{2}(\mathbb{R}^{3})}^{2}\right) ^{3}.
\end{equation*}%
Solving the differential inequality $\frac{d}{dt}f_{n}\leq Cf_{n}^{3}$, with
$f_{n}(t)=\Vert u_{n}(t)\Vert _{H^{2}(\mathbb{R}^{3})}^{2}+\Vert
b_{n}(t)\Vert _{H^{2}(\mathbb{R}^{3})}^{2}$, we get%
\begin{equation*}
f_{n}^{2}(t)\leq \frac{f_{n}^{2}(0)}{1-2Ctf_{n}^{2}(0)}.
\end{equation*}%
Then, using the boundedness of $(f_{n}(0))_{n\in \mathbb{N}}$ by $%
f_{0}=\Vert u_{0}\Vert _{H^{2}(\mathbb{R}^{3})}^{2}+\Vert b_{0}\Vert _{H^{2}(%
\mathbb{R}^{3})}^{2}$ and fixing $0<T^{\ast }<\frac{1}{2Cf_{0}^{2}}$, we
obtain that
\begin{equation}
(u_{n})_{n\in \mathbb{N}}\mbox{ and }(b_{n})_{n\in \mathbb{N}}%
\mbox{
are bounded in }L^{\infty }((0,T^{\ast }),H^{2}(\mathbb{R}^{3})).
\label{af1}
\end{equation}%
Using again (\ref{cach2n}) and equations in (\ref{10z}), we can prove that
\begin{equation}
(\mathcal{J}_{n}u_{n})_{n\in \mathbb{N}}\mbox{ and }(\mathcal{J}%
_{n}b_{n})_{n\in \mathbb{N}}\mbox{
are bounded in }L^{2}((0,T^{\ast }),H^{3}(\mathbb{R}^{3})),  \label{af1-2}
\end{equation}%
\begin{equation}
(\partial _{t}u_{n})_{n\in \mathbb{N}}\mbox{ and }(\partial _{t}b_{n})_{n\in
\mathbb{N}}\mbox{ are bounded in }L^{\infty }((0,T^{\ast }),L^{2}(\mathbb{R}%
^{3})).  \label{af2}
\end{equation}

By (\ref{af1})-(\ref{af2}), there exist a sub-sequence of $%
(u_{n},b_{n})_{n\in \mathbb{N}}$ (still indexed by $n$) and functions $%
u,b\in L^{\infty }((0,T^{\ast }),H^{2}(\mathbb{R}^{3}))\cap L^{2}((0,T^{\ast
}),H^{3}(\mathbb{R}^{3}))$ such that (see \cite{Temam})%
\begin{eqnarray*}
&&(u_{n},b_{n})\overset{n\rightarrow \infty }{\longrightarrow }(u,b)%
\mbox{
weak-* in }L^{\infty }((0,T^{\ast }),H_{\sigma }^{2}(\mathbb{R}^{3})\times
H^{2}(\mathbb{R}^{3})), \\
&&(\mathcal{J}_{n}u_{n},\mathcal{J}_{n}b_{n})\overset{n\rightarrow \infty }{%
\longrightarrow }(u,b)\mbox{
weak in }L^{2}((0,T^{\ast }),H_{\sigma }^{3}(\mathbb{R}^{3})\times H^{3}(%
\mathbb{R}^{3})), \\
&&(u_{n},b_{n})\overset{n\rightarrow \infty }{\longrightarrow }(u,b)%
\mbox{
strong in }L^{2}((0,T^{\ast }),L_{loc}^{2}(\mathbb{R}^{3})\times L_{loc}^{2}(%
\mathbb{R}^{3})), \\
&&(\partial _{t}u_{n},\partial _{t}b_{n})\overset{n\rightarrow \infty }{%
\longrightarrow }(u,b)\mbox{ weak-* in }L^{\infty }((0,T^{\ast }),L^{2}(%
\mathbb{R}^{3})\times L^{2}(\mathbb{R}^{3})).
\end{eqnarray*}%
With these properties, we can apply the weak limit in (\ref{10z}) and prove
that $(u,b)$ is a local strong solution of (\ref{1}) with initial data $%
(u_{0},b_{0})$.

Let us to prove the uniqueness. Suppose that $(u,b)$ and $(v,h)$ are two
weak solutions of (\ref{1}) with the same initial data. Let $U=v-u$ and $%
B=h-b$. So, by inequality (\ref{in5}) (see Remark \ref{Rem-4-10}), we have
that
\begin{equation}
\frac{1}{2}\frac{d}{dt}\left( \Vert U(t)\Vert _{L^{2}}^{2}+\Vert B(t)\Vert
_{L^{2}}^{2}\right) \leq C_{2}\left( \Vert \nabla u(t)\Vert
_{L^{2}}^{4}+\Vert \nabla b(t)\Vert _{L^{2}}^{4}+\Vert \Delta b(t)\Vert
_{L^{2}}^{4}\right) (\Vert U(t)\Vert _{L^{2}}^{2}+\left. \Vert B(t)\Vert
_{L^{2}}^{2}\right) .  \label{aux-uniq-1}
\end{equation}%
Now the uniqueness of solutions in $L^{4}((0,T^{\ast }),H_{\sigma }^{1}(%
\mathbb{R}^{3})\times H^{2}(\mathbb{R}^{3}))$ follows from (\ref{aux-uniq-1}%
) and Gronwall inequality.

In what follows, we prove the blow-up criterion (\ref{asd}). If fact,
suppose that $(u,b)$ satisfies
\begin{equation*}
M:=\int_{0}^{T}\left( \Vert \nabla u(s)\Vert _{L^{2}(\mathbb{R}%
^{3})}^{4}+\Vert \nabla b(s)\Vert _{L^{2}(\mathbb{R}^{3})}^{4}+\Vert \Delta
b(s)\Vert _{L^{2}(\mathbb{R}^{3})}^{4}\right) \,ds<\infty .
\end{equation*}%
Proceeding in an analogous way to the proof of Lemma \ref{9}, we obtain
\begin{equation*}
\Vert u(t)\Vert _{H^{2}(\mathbb{R}^{3})}^{2}+\Vert b(t)\Vert _{H^{2}(\mathbb{%
R}^{3})}^{2}\leq e^{2CM}\left( \Vert u(0)\Vert _{H^{2}(\mathbb{R}%
^{3})}^{2}+\Vert b(0)\Vert _{H^{2}(\mathbb{R}^{3})}^{2}\right) .
\end{equation*}%
So, by using the usual blow-up criterion of time-continuous $H^{2}$%
-solutions, we have that the solution can be extended beyond $T$ (see \cite%
{DC2} and \cite{Temam}).

Finally, let us prove global solutions for small initial data. By inequality
(\ref{cach}) and equivalences (\ref{aux-equiv-1})-(\ref{aux-equiv-2}), there
exist $C_{5}>0$ and $C_{6}>0$ such that
\begin{eqnarray}
\frac{d}{dt}\left( \Vert b(t)\Vert _{H^{2}}^{2}+\Vert u(t)\Vert
_{H^{2}}^{2}\right) &+&C_{5}\left( \Vert u(t)\Vert _{H^{2}}^{2}+\Vert
b(t)\Vert _{H^{2}}^{2}\right)  \notag \\
&\leq &C_{6}\left( \Vert b(t)\Vert _{H^{2}}^{2}+\Vert u(t)\Vert
_{H^{2}}^{2}\right) ^{3}+C_{5}\left( \Vert u(t)\Vert _{L^{2}}^{2}+\Vert
b(t)\Vert _{L^{2}}^{2}\right) .  \label{in20}
\end{eqnarray}%
Suppose that the initial data is small enough to satisfy
\begin{equation*}
\Vert u(0)\Vert _{H^{2}}^{2}+\Vert b(0)\Vert _{H^{2}}^{2}\leq \frac{1}{12}%
\sqrt{\frac{C_{5}}{C_{6}}}.
\end{equation*}%
Let $T^{\ast }>0$ be the supremum over all finite $\widetilde{T}>0$ such
that
\begin{equation*}
\Vert u(t)\Vert _{H^{2}}^{2}+\Vert b(t)\Vert _{H^{2}}^{2}\leq \sqrt{\frac{%
C_{5}}{2C_{6}}},\,\,\forall \,0\leq t\leq \widetilde{T}.
\end{equation*}%
By contradiction, let us assume that $0<T^{\ast }<\infty .$ By (\ref{in20})
we get
\begin{equation*}
\frac{d}{dt}\left( \Vert b(t)\Vert _{H^{2}}^{2}+\Vert u(t)\Vert
_{H^{2}}^{2}\right) +\frac{C_{5}}{2}\left( \Vert u(t)\Vert
_{H^{2}}^{2}+\Vert b(t)\Vert _{H^{2}}^{2}\right) \leq C_{5}\left( \Vert
u(t)\Vert _{L^{2}}^{2}+\Vert b(t)\Vert _{L^{2}}^{2}\right) ,\,\,\forall
\,0\leq t\leq \widetilde{T},
\end{equation*}%
for all $0<$ $\widetilde{T}<T^{\ast }$. Then, Gronwall type inequality and
the time-uniform boundedness $\Vert u(t)\Vert _{L^{2}}^{2}+\Vert b(t)\Vert
_{L^{2}}^{2}\leq \Vert u(0)\Vert _{L^{2}}^{2}+\Vert b(0)\Vert _{L^{2}}^{2}$
give
\begin{eqnarray}
\Vert b(t)\Vert _{H^{2}}^{2}+\Vert u(t)\Vert _{H^{2}}^{2} &\leq &e^{-\frac{%
C_{5}t}{2}}\left( \Vert u(0)\Vert _{H^{2}}^{2}+\Vert b(0)\Vert
_{H^{2}}^{2}+\int_{0}^{t}e^{\frac{C_{5}s}{2}}C_{5}(\Vert u(t)\Vert
_{L^{2}}^{2}+\Vert b(t)\Vert _{L^{2}}^{2})\,ds\right)  \notag \\
&\leq &3(\Vert u(0)\Vert _{H^{2}}^{2}+\Vert b(0)\Vert _{H^{2}}^{2})\leq
\frac{1}{4}\sqrt{\frac{C_{5}}{C_{6}}},\,\,\forall \,0\leq t\leq \widetilde{T}%
,  \label{aux-max-1}
\end{eqnarray}%
for all $0<$ $\widetilde{T}<T^{\ast }.$ In view of the time-continuity of $%
(u,b)$ (see Remark \ref{Rem3}) and $\frac{1}{4}\sqrt{\frac{C_{5}}{C_{6}}}<%
\sqrt{\frac{C_{5}}{2C_{6}}},$ the estimate (\ref{aux-max-1}) contradicts the
maximality of $T^{\ast }.$ So, $T^{\ast }=\infty $ and the solution is
global in time.

\begin{flushright}
  \rule{2mm}{2mm}
 \end{flushright}

%%%%%%%%%%%%%%%%%%%%%%%%%%%%%%%%%%%%%%%%%%%%%%%%%%%%%%%%%%%%%%%%%%%%%%%%%%%%%%%%%%%%%%%%%%%%%%%%%%%%%%%%%

\subsection{\protect\bigskip Proof of Theorem \protect\ref{theo2}}

\hspace{0.6cm}Let $U=v-u$, $B=h-b$ and
\begin{equation*}
L(t)=\Vert \nabla u\Vert _{L^{2}}^{4}+\Vert \nabla b\Vert _{L^{2}}^{4}+\Vert
\Delta u\Vert _{L^{2}}^{4}+\Vert \Delta b\Vert _{L^{2}}^{4}+\Vert \nabla
\times \Delta b\Vert _{L^{2}}^{2}+\Vert div\,\Delta b\Vert _{L^{2}}^{2}.
\end{equation*}%
By Remark \ref{12}, there are $C_{13}>0$, $C_{14}>0$ and $C_{15}>0$ such
that
\begin{eqnarray}
\frac{d}{dt}\left( \Vert U\Vert _{H^{2}}^{2}\right. &+&\left. \Vert B\Vert
_{H^{2}}^{2}\right) +C_{13}\left( \Vert U\Vert _{H^{2}}^{2}+\Vert B\Vert
_{H^{2}}^{2}\right) \leq C_{14}\left( \Vert U\Vert _{H^{2}}^{2}+\Vert B\Vert
_{H^{2}}^{2}\right) ^{3}  \notag \\
&+&C_{15}\left( \Vert U\Vert _{H^{2}}^{2}+\Vert B\Vert _{H^{2}}^{2}\right)
L(t)+C_{13}\left( \Vert U\Vert _{L^{2}}^{2}+\Vert B\Vert _{L^{2}}^{2}\right)
.  \label{aux-proof-teo2}
\end{eqnarray}%
On the other side, using Lemma \ref{9}, we have that if (\ref{444}) holds,
then
\begin{equation*}
\int_{0}^{\infty }L(s)ds<\infty .
\end{equation*}%
Suppose that the initial data $(v(0),h(0))$ is close to $(u(0),b(0))$ to
satisfy (here $C_{2}>0$ is given in (\ref{in5}))
\begin{equation}
\Vert U(0)\Vert _{H^{2}}^{2}+\Vert B(0)\Vert _{H^{2}}^{2}\leq \frac{1}{12}%
\sqrt{\frac{C_{13}}{C_{14}}}\frac{1}{e^{(C_{15}+2C_{2})\int_{0}^{\infty
}L(s)ds}}.  \label{wq1}
\end{equation}%
Let $T^{\ast }>0$ be the supremum over all finite $\widetilde{T}>0$ such
that
\begin{equation*}
\Vert U(t)\Vert _{H^{2}}^{2}+\Vert B(t)\Vert _{H^{2}}^{2}\leq \sqrt{\frac{%
C_{13}}{2C_{14}}},\,\,\forall \,0\leq t\leq \widetilde{T}.
\end{equation*}%
Assume by contradiction that $T^{\ast }<\infty .$ So, for all $0<$ $%
\widetilde{T}<T^{\ast },$ we get%
\begin{eqnarray*}
\frac{d}{dt}\left( \Vert U\Vert _{H^{2}}^{2}\right. +\left. \Vert B\Vert
_{H^{2}}^{2}\right) +\frac{C_{13}}{2}\left( \Vert U\Vert _{H^{2}}^{2}+\Vert
B\Vert _{H^{2}}^{2}\right) &\leq &C_{15}\left( \Vert U\Vert
_{H^{2}}^{2}+\Vert B\Vert _{H^{2}}^{2}\right) L(t) \\
+C_{13}\left( \Vert U\Vert _{L^{2}}^{2}+\Vert B\Vert _{L^{2}}^{2}\right)
\,\,\forall \,0 &\leq &t\leq \widetilde{T}.
\end{eqnarray*}%
By Gronwall type inequality, for all $0<$ $\widetilde{T}<T^{\ast },$ we
obtain
\begin{eqnarray}
\Vert U(t)\Vert _{H^{2}}^{2}+\Vert B(t)\Vert _{H^{2}}^{2} &\leq &e^{-\frac{%
C_{13}}{2}t+C_{15}\int_{0}^{t}L(s)ds}\left( \Vert U(0)\Vert
_{H^{2}}^{2}+\Vert B(0)\Vert _{H^{2}}^{2}\right)  \notag \\
&+&e^{-\frac{C_{13}}{2}t+C_{15}\int_{0}^{t}L(s)ds}C_{13}\int_{0}^{t}e^{\frac{%
C_{13}}{2}s}\left( \Vert U(s)\Vert _{L^{2}}^{2}+\Vert B(s)\Vert
_{L^{2}}^{2}\right) ds  \notag \\
&\leq &\frac{1}{12}\sqrt{\frac{C_{13}}{C_{14}}}+2e^{C_{15}\int_{0}^{\infty
}L(s)ds}\sup_{s>0}\left\{ \Vert U(s)\Vert _{L^{2}}^{2}+\Vert B(s)\Vert
_{L^{2}}^{2}\right\} ,\,\,\forall \,0\leq t\leq \widetilde{T}.  \notag
\end{eqnarray}%
Finally, it follows from (\ref{in5}) and Gronwall inequality that
\begin{equation*}
\sup_{s>0}\left\{ \Vert U(s)\Vert _{L^{2}}^{2}+\Vert B(s)\Vert
_{L^{2}}^{2}\right\} \leq e^{2C_{2}\int_{0}^{\infty }L(s)ds}\left( \Vert
U(0)\Vert _{L^{2}}^{2}+\Vert B(0)\Vert _{L^{2}}^{2}\right) ,\text{ }\forall
\,0\leq t\leq \widetilde{T}.
\end{equation*}%
So
\begin{equation}
\Vert U(t)\Vert _{H^{2}}^{2}+\Vert B(t)\Vert _{H^{2}}^{2}\leq \frac{1}{4}%
\sqrt{\frac{C_{13}}{C_{14}}},\,\,\forall \,0\leq t\leq \widetilde{T},
\label{wq2}
\end{equation}%
for all $0<$ $\widetilde{T}<T^{\ast }.$ The estimate (\ref{wq2}) contradicts
the maximality of $T^{\ast }$ because $\frac{1}{4}\sqrt{\frac{C_{13}}{C_{14}}%
}<\sqrt{\frac{C_{13}}{2C_{14}}}$ and $(U,B)$ is time-continuous (see Remark %
\ref{Rem3}). It follows that $T^{\ast }=\infty .$ As $(u,b)$ is a global
solution in $H^{2}(\mathbb{R}^{3})$, the above inequality implies that $%
(v,h) $ is as well. Furthermore, by repeating steps between (\ref{wq1})-(\ref%
{wq2}), one can check that for initial data less than $\delta $, where $%
0<\delta <\frac{1}{12}\sqrt{\frac{C_{13}}{C_{14}}}\frac{1}{%
e^{(C_{15}+2C_{2})\int_{0}^{\infty }L(s)ds}}$, we can take $M(\delta
)=3\delta e^{(C_{15}+2C_{2})\int_{0}^{\infty }L(s)ds}$, for $M(\delta )$ as
in the statement of the theorem. This concludes the proof.

\begin{flushright}
  \rule{2mm}{2mm}
 \end{flushright}

%%%%%%%%%%%%%%%%%%%%%%%%%%%%%%%%%%%%%%%%%%%%%%%%%%%%%%%%%%%%%%%%%%%%%%%%%%%%%%%%%%%%%%%%%%%%%%%%%%%%%%%%%

%%%%%%%%%%%%%%%%%%%%%%%%%%%%%%%%%%%%%%%%%%%%%%%%%%%%%%%%%%%%%%%%%%%%%%%%%%%%%%%%%%%%%%%%%%%%%%%%%%%%%%%%%
%%%%%%%%%%%%%%%%%%%%%%%%%%%%%%%%%%%%%%%%%%%%%%%%%%%%%%%%%%%%%%%%%%%%%%%%%%%%%%%%%%%%%%%%%%%%%%%%%%%%%%%%%

\end{document}